\documentclass[a4paper]{article}
\usepackage[utf8]{inputenc}
\usepackage{amsmath,amsthm,amssymb}
\usepackage{enumerate}
\usepackage{float}
\usepackage{cite} 
\usepackage{graphicx,xcolor,cancel}
\usepackage[toc,page]{appendix}
\usepackage{comment}
\usepackage{mathtools}
\usepackage[bb=boondox,bbscaled=1.2,cal=boondoxo]{mathalfa}
\usepackage{braket,mathdots}
\usepackage{array}
\usepackage{relsize}

\usepackage[bookmarks=false,hidelinks]{hyperref} 
\usepackage[top=2.8cm,bottom=2.8cm,left=2.6cm,right=2.6cm]{geometry}
\usepackage{multicol}
\usepackage{tikz}
\usetikzlibrary{shapes.geometric,decorations.markings,intersections,calc,positioning,arrows.meta}

\theoremstyle{definition}
\newtheorem{remark}{Remark}[section]
\newtheorem{theorem}{Theorem}[section]
\newtheorem{corollary}{Corollary}[section]
\newtheorem{prop}{Proposition}[section]

\numberwithin{equation}{section}

\expandafter\def\expandafter\normalsize\expandafter{%
    \normalsize%
    \setlength\abovedisplayskip{6pt}%
    \setlength\belowdisplayskip{6pt}%
    \setlength\abovedisplayshortskip{4pt}%
    \setlength\belowdisplayshortskip{4pt}%
}

\title{Skew-orthogonal polynomials for a quartic Freud weight: \\ 
two classes of quasi-orthogonal polynomials}
\author{Costanza Benassi\footnote{Corresponding author.}\\[-.5ex]
\small{\url{costanza.benassi@northumbria.ac.uk}} \\[-.5ex]
{\small ORCID iD: \url{https://orcid.org/0000-0001-7254-8559}}\\[.5ex]
{ \small School of Engineering, Physics and Mathematics, Northumbria University,} {\small Newcastle upon Tyne, UK}
\\
Marta Dell'Atti \\[-.5ex]\small{\url{m.dell-atti@uw.edu.pl}} \\[-.5ex]
{\small ORCID iD: \url{https://orcid.org/0000-0001-6484-1739}} \\[.5ex]
{ \small University of Warsaw, Institute of Mathematics, Banacha 2, Warsaw, Poland}}
\date{}

\allowdisplaybreaks

\begin{document}

\maketitle
\begin{abstract}
 This work is a thorough investigation of skew-orthogonal polynomials with respect to a quartic Freud weight.  We provide an explicit method to evaluate skew-orthogonal polynomials of any degree as linear combinations of orthogonal polynomials.  The coefficients of these combinations can be evaluated via novel recursive relations. Moreover, we observe that skew-orthogonal polynomials with even and odd degree constitute two families of quasi-orthogonal polynomials with respect to two different semi-classical Laguerre weights,  and we provide the first instance of closed recursive relations involving skew-orthogonal polynomials only. 
\end{abstract}

\section{Introduction}

Orthogonal polynomials play a pivotal role in a number of fields in mathematics, including analysis and mathematical physics, and their properties are deeply related to those of the weight defining the relevant scalar product~\cite{szeg1939orthogonal,szego1975orthogonal}. In this work, we focus on a quartic Freud weight, i.e.\ the exponential of a quartic potential of the form $w(x;t) = \exp(-x^4 + t x^2)$. Polynomials orthogonal with respect to this weight and its generalisations to higher degree potentials have been object of much research effort 
~\cite{FilipukVanAsscheZhang,Clarkson_13,ClaJorKel2016, ClarJordaLour2023, Clarkson25}. 

The aim of this work is to investigate the relationship between Freud orthogonal polynomials and their skew-orthogonal counterparts. 
Indeed, just as  for a given weight it is possible to define a (symmetric) scalar product on an appropriate space of functions {\color{black}and related orthogonal polynomials}, similarly one can define a skew-symmetric scalar product, inducing a family of skew-orthogonal polynomials. This has emerged naturally in the context of random matrix theory, and skew-orthogonal polynomials play a crucial role in the connection between Orthogonal and Symplectic Ensembles of random matrices and integrable lattices~\cite{Dyson72, BrezinSkew, ghosh, Mehta3, Mehta4, Meh2004, NagaoII, NagaoIII}. They have been further investigated in this context, and their properties {\color{black}used} to explore the related random matrix models \cite{vanM_Pfaff_skew, Adler2000_skew, vanMoerbekenotes, vanM_Pfaff_tau_function, adler2002, BDM2}. Of particular relevance for our purposes is \cite{Adler2000_skew}, where a mapping  between skew-orthogonal and orthogonal polynomials was established in a general framework. In the present paper, we build upon these results and characterise this mapping in the case of a quartic Freud weight, expressing {\color{black}polynomials skew-orthogonal} with respect to $w(x;t)$ as finite linear combinations of polynomials orthogonal with respect to the weight $w(x;t)^2$. We show that each skew-orthogonal polynomial of odd degree is given as a combination of three odd orthogonal ones, whilst skew-orthogonal polynomials of even degree are expressed in terms of two even orthogonal ones. The coefficients of these linear combinations can be evaluated explicitly as they are shown to satisfy recursive relations, allowing us to build skew orthogonal polynomials for any value of the parameter $t$ in the weight. It is worth mentioning that the skew-symmetric scalar product that we use to generate the skew-orthogonal polynomials exhibits the so-called biorthogonality property with a skew-symmetric kernel. Biorthogonal pairs of sequences were already mentioned in the classical work by Szeg\"o~\cite{szeg1939orthogonal,szego1975orthogonal}, whereas in~\cite{norset} the notion of a sequence of polynomials biorthogonal with respect to different {\color{black}scalar products} was introduced. More recently, in the context of random matrix ensembles, biorthogonal polynomials were introduced in~\cite{BertolaBiorto1,BertolaBiorto2} for products built with totally positive kernels, and extended to skew-symmetric kernels in~\cite{LiShenXianYu}. 

Another class of weights deeply connected to our {\color{black}setting} is semi-classical Laguerre weights. They constitute a family of weights of the form $w_{\lambda}(z;t) = z^\lambda \exp(-2z^2 + 2 tz)$ depending on a parameter $\lambda > -1$. It is well known that Freud orthogonal polynomials can be mapped into semi-classical Laguerre orthogonal ones -- in particular, {\color{black}even-degree} ones constitute a family of polynomials orthogonal with respect to $w_{-\frac{1}{2}}(z;t)$ {\color{black} with $z=x^2$}, whilst {\color{black}odd-degree} ones {\color{black}are orthogonal} with respect to $w_{\frac{1}{2}}(z;t)$ (see for example \cite{szeg1939orthogonal, Carlitz}, or more recently \cite{FilipukVanAsscheZhang}). In this work we prove that Freud {\color{black}even-degree} skew-orthogonal polynomials correspond to a linear combination of two subsequent {\color{black}polynomials orthogonal} with respect to {\color{black}Laguerre} $w_{-\frac{1}{2}}(z;t)$, whilst {\color{black}odd-degree} ones can be mapped {\color{black}onto} a linear combination of three subsequent {\color{black}polynomials} orthogonal with respect to {\color{black}Laguerre} $w_{\frac{1}{2}}(z;t)$.

These results put skew-orthogonal polynomials in the framework of so-called quasi-orthogonal polynomials, i.e.\ polynomials defined as appropriate linear combinations of orthogonal ones. The first examples of quasi-orthogonal polynomials are found in \cite{riesz,fejer}, while a more general definition is given in \cite{chihara}. Recent work on quasi-orthogonal polynomials has focused on several aspects: the conditions for orthogonality (or obstructions to it)~\cite{Draux01092016,IsmWan,JordaanJost}, the derivation of recurrence relations~\cite{Bracciali}, the analysis of their zeros~\cite{DriverJordaan}, and their connections to Painlev\'e equations~\cite{FILIPUKquasiortho}. Our results prove that {\color{black}polynomials skew-orthogonal} with respect to the Freud weight $w(x;t)$ constitute two families of {\color{black}polynomials quasi-orthogonal} with respect to two different semi-classical Laguerre weights $w_\lambda(z;t)$. This allows us to use the machinery and results developed in the context of quasi-orthogonal polynomials, and apply it to Freud skew-orthogonal ones. In particular, it was shown in~\cite{IsmWan} that, whilst {\color{black}polynomials orthogonal} with respect to a given weight fulfil a closed three-point relation with constant coefficients except one linear in the independent variable (see e.g.\ \cite{BonNev}), quasi-orthogonal polynomials satisfy three-point relations with all coefficients depending in a non-trivial way on the independent variable, including higher degree terms {\color{black}(nonlinear)}. {\color{black}In this work, we} exploit the quasi-orthogonality of Freud skew-orthogonal polynomials to provide new explicit closed {\color{black}three-{\color{black}point}} expressions. 

This paper is structured as follows. In Section \ref{sec:Freud} we collect some previously known results concerning orthogonal, quasi-orthogonal and skew-orthogonal polynomials, with focus on the quartic Freud weight. 
Sections \ref{sec:transition} and \ref{sec:skew&quasi} include the novel results of this work. Section \ref{sec:transition} {\color{black}focuses} on the mapping between skew-orthogonal and orthogonal polynomials in the specific case of Freud weights. We provide an explicit formulation of this mapping, proving that it leads to skew-orthogonal polynomials being expressed as suitable linear combinations of orthogonal ones, and hence {\color{black}establishing} their quasi-orthogonality. The coefficients of these combinations satisfy novel recursive relations which we can solve iteratively with appropriate software, and we provide relevant plots to exemplify their behaviour (see Fig.~\ref{fig:plots}). The results {\color{black} just mentioned} are formulated in Theorem \ref{thm:skew_quasi} and Corollary \ref{cor:quasi_Laguerre}. 
In Section \ref{sec:skew&quasi} we exploit the quasi-orthogonality of skew-orthogonal polynomials of even and odd degree and, following the approach proposed in \cite{IsmWan}, we provide two new explicit closed expressions for skew-orthogonal polynomials, one for those of even degree and another one for those of odd degree. This is specified in Theorem \ref{thm:even_rec} and Corollary \ref{cor:rec_Laguerre_even} for even degree polynomials, and in Theorem \ref{thm:odd_rec} and Corollary \ref{cor:rec_Laguerre_odd} for odd degree ones. {\color{black} Appendix \ref{app:initial} includes some calculations omitted from Section \ref{sec:transition}.  Explicit expressions for the first few skew-orthogonal polynomials are provided in Appendix \ref{app:formskew}.} Finally, in Appendix~\ref{sec:skewHermite} we look at a simpler set-up and provide closed recursive expressions for the skew-orthogonal version of the Hermite polynomials.

\section{Quartic Freud weight and related polynomials}
\label{sec:Freud}

In this section we collect some results regarding orthogonal, quasi-orthogonal and skew-orthogonal polynomials, with a particular focus on polynomials related to the quartic Freud weight~\cite{Freud1976}
\begin{equation}
\label{eq:freud_weight}
    w(x;t) = \exp(-V(x;t))= \exp(-x^4 +t x^2). 
\end{equation}
In Section \ref{sec:orthointro} we recall some well established features of orthogonal polynomials, including relevant recursive relations and the relationship between Freud and semi-classical Laguerre polynomials. 
Section \ref{sec:quasiorthointro} is devoted to a brief introduction to quasi-orthogonal polynomials and some key results concerning their underlying recursive relations~\cite{chihara, IsmWan}.
In Section \ref{sec:skewintro} we introduce skew-orthogonal polynomials and the general framework for the mapping between skew-orthogonal and orthogonal ones, as detailed in \cite{Adler2000_skew}, which underpins further developments of this work.

\subsection{Orthogonal polynomials}
\label{sec:orthointro}

Given a generic weight $\omega(x):\mathbb{R}\rightarrow \mathbb{R}^+$, we define the symmetric scalar product $(~\cdot~,~\cdot~)$ {\color{black}{\color{black} of} suitable functions $f$ and $g$} as
\begin{equation}
\left(\,f,\, g \,\right) = \int_{\mathbb{R}} f(x)\,  g(x) \, \omega(x) \, {\rm d}x .
    \label{eq:sym_prod_generic}
\end{equation}
{\color{black}{\color{black}Within }this work, we assume $f$ and $g$ to be polynomials in $x$}.
We consider the semi-infinite vector of monic polynomials $P(x)=(P_0(x), P_1(x),\dots)^{\top}$ with elements
\begin{equation}
P_{j}(x) = x^{j}+a_{j-1}\,x^{j-1}+\dots+a_{0}, \qquad a_{j} \equiv 1, ~ a_{j-1}, \dots, a_{0} \in \mathbb{R}, ~ \forall ~j \in \mathbb{N}_0, 
\end{equation}
which are orthogonal with respect to the weight $\omega(x)$, i.e.\ such that
\begin{equation}
    \left(\,P_n,\, P_m\,\right) = \int_{\mathbb{R}} P_n(x)\, P_m(x) \,\omega(x) \, {\rm d}x = h_n \, \delta_{n,m}, \qquad \forall ~ m,n\in\mathbb{N}_0,
    \label{eq:monicgeneric}
\end{equation}
where $\{h_j\}_{j\in\mathbb{N}_0}$ are appropriate constants depending on $\omega(x)$. The orthonormal counterpart of $P(x)$ is the vector $p(x)=(p_0(x),p_1(x), \dots)^{\top}$ with elements satisfying
\begin{equation}
    (\,p_n,\,p_m\,) = \int_{\mathbb{R}}p_n(x)\,p_m(x)\,\omega(x)\,{\rm d}x = \delta_{n,m}, \qquad \forall ~ m,n\in\mathbb{N}_0.
\end{equation}
The orthonormal polynomials are mapped into monic ones via $P_n(x)=\sqrt{h_n}\,p_n(x)$ for any $n\in\mathbb{N}_0$.

A defining characteristic of orthogonal polynomials is that they satisfy the three-term recurrence relation 
\begin{equation}
    x\,P_{n}(x) = P_{n+1}(x) + \alpha_n\,P_n(x)+\beta_n\,P_{n-1}(x), \qquad n \in \mathbb{N}_0,
    \label{eq:recgeneral}
\end{equation}
(see e.g.\ \cite{BonNev}) where the recurrence coefficients $\alpha_n$, $\beta_n$ are defined as 
\begin{equation}\label{eq:alpha&beta}
    \alpha_n = \frac{(\,xP_n,P_n\,)}{(\,P_n,P_n\,)}  , \qquad \beta_n =  \frac{(\,xP_n,P_{n-1}\,)}{(\,P_{n-1},P_{n-1}\,)} ,
\end{equation}
in terms of the product \eqref{eq:sym_prod_generic}.
Let us now consider the case where $\omega(x)$ 
{\color{black} is given by $\omega(x) := \omega(x;\,t) = w(x;t)^2$}, with $w(x;t)$ the Freud weight in~\eqref{eq:freud_weight}, i.e.\ 
\begin{equation}\label{eq:omega_freud}
    \omega(x;t) = w(x;t)^2 = \exp(-2x^4+2tx^2).
\end{equation}
The reason behind this choice {\color{black} of weight} will be more transparent in Section~\ref{sec:skewintro}. {\color{black}Note that for this specific weight we have introduced the notation $\omega(x;\,t)$, to stress its dependency on the parameter $t$.} {\color{black}Therefore,} all relevant quantities in this framework depend on the parameter $t$, as this appears in the weight $w(x;t)$. The corresponding symmetric inner product in~\eqref{eq:sym_prod_generic} will be denoted as $(~\cdot~,~\cdot~)_t$. As the Freud weight is even (i.e.\ it is invariant under $x\rightarrow -x$), it follows that $\alpha_n(t) = 0$ for all $n\geq 0$ in~\eqref{eq:alpha&beta}. Therefore the three-term relation~\eqref{eq:recgeneral} involves $\beta_n(t)$ only and takes the form 
\begin{equation}
\label{eq:3pts}
 x\, P_n(x;t) = P_{n+1}(x;t) + \beta_n(t)\, P_{n-1}(x;t),{\color{black}\qquad n\geq 0},\\
\end{equation}
{\color{black} with $P_{-1}(x;\,t) =0$ and $P_0(x;\,t)=1$.} The recurrence coefficients are given in terms of ratios of the normalisation coefficients of the orthogonal polynomials as 
\begin{equation}\label{eq:beta_def}
    \beta_n(t) = \frac{h_{n}(t)}{h_{n-1}(t)}=\frac{\left(P_n, P_n\right)_t}{\left(P_{
   n-1}, P_{n-1}\right)_t
    },\qquad{\color{black} n>0,\qquad \beta_0(t) = 0.} 
\end{equation}
{\color{black} We also point out that $\beta_n(t)$ is well defined for $n>0$ since $h_{n-1}(t)\neq 0$ as per equation \eqref{eq:monicgeneric}.} From the classical results in~\cite{BonNev}, the so-called structure relation is   
\begin{equation}
\label{eq:derivative}    
\begin{split}
    P_{n}'(x;t) &= n\, P_{n-1}(x;t) + 8\frac{h_n(t)}{h_{n-3}(t)} P_{n-3}(x;t) \\&=n\, P_{n-1}(x;t) + 8 \,\beta_n(t) \beta_{n-1}(t) \beta_{n-2}(t) P_{n-3}(x;t),
\end{split}
\end{equation}
where $P'_n(x;t)$ denotes the derivative of $P_n(x;t)$ with respect to $x$ (see~\cite{VanAssche_2017}).
The compatibility between the recurrence relation~\eqref{eq:3pts} and the structure relation~\eqref{eq:derivative} for the quartic Freud weight gives rise to the recurrence relation for
the coefficients~$\{\beta_j(t)\}_{j\in\mathbb{N}_0}$, i.e.\ 
\begin{equation}
    n = 8\, \beta_n(t) \big(\beta_{n+1}(t)+\beta_n(t) + \beta_{n-1}(t)\!\big)-4\,t\,\beta_n(t)\,{\color{black},}
    \label{eq:recurrence_beta}
\end{equation}
{\color{black} with initial condition $\beta_0(t) = 0$ and $\beta_1(t) = h_1(t)/h_0(t)$.}\footnote{\color{black} At $t=0$ we have 
\begin{equation*}
 \beta_0(0)=0, \qquad   \beta_1(0) = \frac{\Gamma(\frac{3}{4})}{\sqrt{2}\, \Gamma(\frac{1}{4})}, 
\end{equation*}
where we considered the standard integral representation of $\Gamma(\zeta)=\int_0^\infty z^{\zeta-1}\exp(-z)dz$. For a generic time $t$, we can make use of the parabolic cylinder special functions $D_\nu(\zeta)$ as in~\cite[pg.\ 13]{ClaJorKel2016}, with integral representation 
\begin{equation*}
    D_{\nu}(\zeta) = \frac{\exp\bigl(-\frac{\zeta^2}{4}\bigr)}{\Gamma(-\nu)}\int_0^\infty u^{-\nu-1} \exp\bigl(-\frac{1}{2}u^2-\zeta u\bigr) \text{d}u. 
\end{equation*}
With $u=\sqrt{2c}v$ and $\zeta = - \gamma/\sqrt{2c}$, we can write the integral expressions for {\color{black} both} $h_1(t)$ and $h_0(t)$ in terms of these special functions, via the change of variables $x^2 = z$ and identifying $c = 2$, $\gamma = 2t$ and the appropriate $\nu$. The starting points for the recurrence relation~\eqref{eq:recurrence_beta} are then 
\begin{equation*}
    \beta_0(t) = 0, \qquad \beta_1(t) = \frac{D_{-3/2}(-t)}{4 D_{-1/2}(-t)}. 
\end{equation*}
}
Incidentally, in~\cite{FokasIts} the authors named this equation discrete Painlev\'e I (d-$\text{P}_{\text{I}}$). The recurrence relation~\eqref{eq:recurrence_beta} will play a crucial role in the following sections.

We note here that, throughout this work, most quantities will depend on the variables $x$ or $t$, or both. When it is clear from the context, we will often drop this explicit dependency to enhance legibility.

The Freud weight~\eqref{eq:freud_weight} {\color{black}-- or equivalently its re-scaled version in~\eqref{eq:omega_freud} --} is the simplest example of semi-classical weight on the real line. 
A weight~$\omega(x)$ in~\eqref{eq:sym_prod_generic} is said classical if it satisfies the Pearson equation (see for example~\cite{Clarkson_13, VanAssche_2017} and references therein)
\begin{equation}
\label{eq:pearson}
    \frac{\text{d}}{\text{d}x} \big( \sigma (x)\,\omega(x) \big) = \tau(x) \,\omega(x),
\end{equation}
where $\sigma(x)$ is a polynomial with degree $\text{deg}(\sigma)\le 2$, and $\tau(x)$ a polynomial with $\text{deg}(\tau)=1$. Classes of polynomials orthogonal with respect to the weight are then called classical and notable examples are Hermite, Laguerre and Jacobi orthogonal polynomials. On the other hand, so-called semi-classical weights satisfy the Pearson equation~\eqref{eq:pearson} where either $\sigma(x)$ is a polynomial with $\text{deg}(\sigma)>2$ or $\tau(x)$ a polynomial with $\text{deg}(\tau)\neq 1$. In the case of our weight $\omega(x;t)$ in~\eqref{eq:omega_freud} 
\begin{equation}
    \sigma(x) = 1 , \qquad \tau(x) = -4x^3 +2tx. 
\end{equation}
As already noticed in~\cite{FilipukVanAsscheZhang} in terms of recurrence coefficients, and generalising results from \cite{szeg1939orthogonal, Carlitz}, we remark that the weight $\omega(x;t)$ in~\eqref{eq:omega_freud} is related to semi-classical Laguerre polynomials orthogonal with respect to the weight
\begin{equation}\label{eq:semi-Laguerre_general}
    w_{\lambda}(z;t) = z^{\lambda} \exp(-2z^2+ 2t z), \qquad z>0,~~ \lambda > -1.
\end{equation}
This relationship is made evident by the change of variable $z:=x^2$.  
Firstly, we remark that it follows straightforwardly from the symmetry of the  weight \eqref{eq:omega_freud} that any even degree polynomial $P_{2n}(x;t)$ includes terms with even powers of the independent variable $x$ only, whilst any odd degree polynomial $P_{2n+1}(x;t)$ comprises terms with odd powers of $x$ only. {\color{black} As a consequence, orthogonal polynomials appearing in~\eqref{eq:3pts} with even and odd index can be mapped onto polynomials which are orthogonal with respect to the weights $\omega(z;t)=w_{-\frac{1}{2}}(z;t)$ and $\omega(z;t)=w_{\frac{1}{2}}(z;t)$ respectively\footnote{\color{black}We remark that the standard form of the semi-classical Laguerre weight is 
\begin{equation*}
\overline{w}_{\lambda
 }(u;t) = u^{\lambda} \, \exp({-u^2+\overline{t}\,u})   .  
 \label{eq:semi-Laguerre_standard}
\end{equation*}
However, this can be found from \eqref{eq:semi-Laguerre_general} with a simple change of variables $z:= \frac{u}{\sqrt{2}}$ and $\overline{t}:=\sqrt{2}\,t$. The semi-classical Laguerre weight in its standard formulation has been extensively studied in  \cite{Clarkson_13}.}.} 

To further clarify this, let us first consider an even degree polynomial $P_{2n}(x;t)$. As it only includes terms with even powers of $x$, we can set $z:= x^2$ 
and define 
\begin{equation}\label{eq:p_ortho_laguerre}
    \widehat{P}_{n}(z;t):=P_{2n}(x;t).
\end{equation}
It is a simple calculation to show that the polynomials $\widehat{P}(z;t) = (\widehat{P}_{0}(z;t),\,\widehat{P}_{1}(z;t),\,\dots)^{\top}$ are monic orthogonal with respect to the semi-classical Laguerre weight \eqref{eq:semi-Laguerre_general}
 with $\lambda=-\frac{1}{2}$, and they satisfy the three-point relation 
\begin{equation}
\label{eq:rec_Phat}
   z\, \widehat{P}_{n}(z;t) = \widehat{P}_{n+1}(z;t) + \widehat{a}_n(t)\,\widehat{P}_{n}(z;t) + \widehat{b}_n(t)\,\widehat{P}_{n-1}(z;t) , 
\end{equation}
with
\begin{equation}\label{eq:laguerre_coefficients}
        \widehat{a}_n(t) := \beta_{2n}(t)+\beta_{2n+1}(t),\qquad\qquad
        \widehat{b}_n(t) := \beta_{2n-1}(t)\,\beta_{2n}(t).
\end{equation}
This follows from applying twice the three-point relationship \eqref{eq:3pts} to a polynomial of even degree $P_{2n}(x;t)$, giving
\begin{equation}\label{eq:two_step_three_terms}
        x^2\,P_{2n}(x;t)=P_{2n+2}(x;t) + (\beta_{2n}(t)+\beta_{2n+1}(t))P_{2n}(x;t) + \beta_{2n-1}(t)\,\beta_{2n}(t)\, P_{2n-2}(x;t) , 
    \end{equation}
and then mapping each even degree polynomial $P_{2j}(x;t)$ into $\widehat{P}_j(z;t)$ as discussed.

We now turn to the case of polynomials with odd degree $P_{2n+1}(x;t)$ orthogonal with respect to the weight~$\omega(x;t)$ in~\eqref{eq:omega_freud}. As before, let $z = x^2$, and define the vector of polynomials $\widetilde{P}(z;t)=(\widetilde{P}_0(z;t),\,\widetilde{P}_1(z;t),\,\dots)^{\top}$ as
\begin{equation}
\label{eq:Ptildedef}
   \widetilde{P}_n(z;t):= \frac{P_{2n+1}(x;t)}{x},
\end{equation}
exploiting the fact that $P_{2n+1}(x;t)$ comprises terms with odd powers of $x$ only.
It is straightforward to check that polynomials $\widetilde{P}(z;t)$ are monic orthogonal with respect to the semi-classical Laguerre weight~\eqref{eq:semi-Laguerre_general} with $\lambda = \frac{1}{2}$. They satisfy the three-point relation 
\begin{equation}
\label{eq:rec_Ptilde}
   z\, \widetilde{P}_{n}(z;t) = \widetilde{P}_{n+1}(z;t) + \widetilde{a}_n(t)\,\widetilde{P}_{n}(z;t) + \widetilde{b}_n(t)\,\widetilde{P}_{n-1}(z;t) . 
\end{equation}
with
\begin{equation}\label{eq:laguerre_coefficientsodd}
    \begin{split}
        \widetilde{a}_n(t) &:= \beta_{2n+2}(t)+\beta_{2n+1}(t), \qquad \qquad 
        \widetilde{b}_n(t) = \beta_{2n+1}(t)\,\beta_{2n}(t).
     \end{split}
\end{equation}
Equation \eqref{eq:rec_Ptilde} follows from applying \eqref{eq:3pts} twice to an odd degree polynomial $P_{2n+1}(x;t)$, finding
\begin{equation}
    x^2 P_{2n+1}(x;t) = P_{2n+3}(x;t) + (\beta_{2n+2}(t) + \beta_{2n+1}(t))P_{2n+1}(x;t) + \beta_{2n+1}(t)\beta_{2n}(t)P_{2n-1}(x;t) 
\label{eq:x2odd}
\end{equation}
and then mapping each odd degree polynomial $P_{2j+1}(x;t)$ into $\widetilde{P}_{j}(z;t)$ according to~\eqref{eq:Ptildedef}.

{\color{black}
\begin{remark} 
We note that it is possible to invert  expressions~\eqref{eq:laguerre_coefficients} to formulate $\{\beta_{j}(t)\}_{j\in\mathbb{N}_0}$ in terms of $\left\{\widehat{a}_j(t)\right\}_{j\in\mathbb{N}_0}$ and $\{\widehat{b}_j(t)\}_{j\in\mathbb{N}_0}$. We find a formulation in continued fractions for the case of even index $\beta_{2n}(t)$, i.e.\ 
\begin{equation}
\begin{split} 
    \beta_{2n} &= \widehat{a}_{n} + \frac{\widehat{b}_{n+1} }{\widehat{a}_{n+1} + \dfrac{\widehat{b}_{n+2} }{\widehat{a}_{n+2} + \dfrac{\widehat{b}_{n+3} }{\widehat{a}_{n+3} + \dfrac{\widehat{b}_{n+4} }{\dots}}}},
\end{split} 
\end{equation}
via the iteration of the following expression, which follows directly from \eqref{eq:laguerre_coefficients} 
\begin{equation}
  \beta_{2n} = \widehat{a}_{n} + \frac{\widehat{b}_{n} }{\beta_{2n+2}}.
\end{equation}
For odd indices $\beta_{2n-1}(t)$ we find an analogous expression, i.e.\ 
\begin{equation}
    \beta_{2n-1} = \frac{\widehat{b}_n}{\widehat{a}_n-\dfrac{\widehat{b}_{n+1}}{\widehat{a}_{n+1}-\dfrac{\widehat{b}_{n+2}}{\widehat{a}_{n+2}-\dfrac{\widehat{b}_{n+3}}{\widehat{a}_{n+3}-\dots}}}},
\end{equation}
from the iteration of 
\begin{equation}
    \beta_{2n-1} = \frac{\widehat{b}_{n} }{\widehat{a}_{n} -\beta_{2n+1} },
\end{equation}
which is a straightforward consequence of \eqref{eq:laguerre_coefficients}.

It can be shown that with a similar approach we can invert~\eqref{eq:laguerre_coefficientsodd} and formulate $\{\beta_{j}(t)\}_{j\in\mathbb{N}_0}$ as different expressions in continued fractions in terms of $\left\{\widetilde{a}_j(t)\right\}_{j\in\mathbb{N}_0}$ and $\{\widetilde{b}_j(t)\}_{j\in\mathbb{N}_0}$.
\end{remark}
}

\subsection{Quasi-orthogonal polynomials} 
\label{sec:quasiorthointro}

In this section we briefly introduce the concept of quasi-orthogonal polynomials~\cite{chihara}.
Let us consider monic polynomials $P(x)$ orthogonal with respect to the generic weight $\omega(x)$ as per equation~\eqref{eq:monicgeneric}.
For any given integers $k\geq 1$ and $r\geq 1$, a family of monic polynomials $R(x) = \left(R_0(x), R_1(x), \dots\right)^\top$ {\color{black}quasi-orthogonal with respect to $\omega(x)$ of order~$(k,r)$ is} defined as
\begin{equation}
    R_n(x) = P_n(x) + \sum_{j=1}^k {\color{black}\mathcal{c}}_{n,j}\,P_{n-jr}(x)\qquad \forall ~n\in\mathbb{N}_0,   \label{eq:quasiortho}
\end{equation}
where ${\color{black}\mathcal{c}}_{n,j}$ are parameters not depending on the independent variable $x$, and with the convention that $P_j(x) = 0$ for $j<0$.
The simplest example is given by quasi-orthogonal polynomials of order $(1,1)$ (also known as first-order), of the form
\begin{equation}
\label{eq:quasiortho1st}
    R_n(x) = P_n(x) + {\color{black}\mathcal{c}}_n \, P_{n-1}(x).
\end{equation}
In general, quasi-orthogonal polynomials $R(x)$ in~\eqref{eq:quasiortho} satisfy
\begin{equation}
    \int R_n(x) \, R_m(x) \, \omega(x) \,\text{d}x = 0, \qquad \text{  for  }m\neq n \pm j r,~j\in\{0,1, \dots,k\}  .   \label{eq:quasiortho_int}
\end{equation}
One question that has been much considered in the literature is under which conditions quasi-orthogonal polynomials are orthogonal with respect to an appropriate weight. This has been explored in~\cite{IsmWan} in the case of polynomials of the form~\eqref{eq:quasiortho1st}, and in a more general framework in~\cite{Bracciali}. 
As previously mentioned, any family of orthogonal polynomials satisfies a three-term recurrence relation of the form~\eqref{eq:recgeneral} -- therefore, the problem of orthogonality for quasi-orthogonal polynomials can be formulated as finding conditions such that they satisfy a suitable three-point relation. In particular, as first pointed out in \cite{chihara} and further specified in \cite{IsmWan}, polynomials of the class  \eqref{eq:quasiortho1st} always satisfy a three-point relation of the form
\begin{equation}
    \ell_{n-1}(x)\, R_{n+1}(x) = \big({\color{black}s}_n(x)\, \ell_{n-1}(x)-{\color{black}\mathcal{c}}_{n-1} \,\ell_n(x) \big)R_n(x) - \beta_{n-1}\,\ell_n(x)\, R_{n-1}(x),
    \label{eq:quasiortho-rec}
\end{equation}
{\color{black}in terms of the auxiliary functions $s_n$ and $\ell_n$ defined as}
\begin{equation}
 {\color{black}s}_n(x):= x - \alpha_n + {\color{black}\mathcal{c}}_{n+1},\qquad \ell_n(x):= {\color{black}\mathcal{c}}_n\, {\color{black}s}_n(x) + \beta_n.    
 \label{eq:r&l}
\end{equation}
Here, $\alpha_n$ and $\beta_n$ are the three-point coefficients defined in \eqref{eq:alpha&beta}, $h_n$ is the normalising constant in \eqref{eq:monicgeneric}, ${\color{black}\mathcal{c}}_n$ appears in \eqref{eq:quasiortho1st}, and conventionally $R_j(x) = 0$ if $j<0$.  It is important to notice that the three-point relation in~\eqref{eq:quasiortho-rec} is rather different from the ones satisfied by orthogonal polynomials~\eqref{eq:recgeneral}, as the coefficients of the recurrence relation depend on $x$ in a non-trivial way (linearly, quadratically and linearly for the coefficients of $R_{n+1}$, $R_n$ and $R_{n-1}$ respectively). In \cite{IsmWan} the authors discuss conditions for~\eqref{eq:quasiortho-rec} to reduce to a form similar to \eqref{eq:recgeneral}, i.e.\ for quasi-orthogonal polynomials \eqref{eq:quasiortho1st} to be orthogonal. These amount to the following recurrence relation and inequality holding\footnote{{\color{black}Here it is understood that $\mathcal c_n\neq0$, as the family is assumed
to be genuinely quasi-orthogonal of order one~\cite{IsmWan}.}}:
\begin{subequations}\label{eq:ortho-conditions}
\begin{align}
    {\color{black}\mathcal{c}}_{n+1} = {\color{black}\mathcal{c}}_n + \frac{\beta_{n-1}}{{\color{black}\mathcal{c}}_{n-1}}-\frac{\beta_n}{{\color{black}\mathcal{c}}_n} + \alpha_n-\alpha_{n-1},\qquad&n\geq 2,\\
   \beta_n + {\color{black}\mathcal{c}}_n\left(\alpha_{n-1}-\alpha_n-{\color{black}\mathcal{c}}_n +{\color{black}\mathcal{c}}_{n+1}\right)>0,\qquad&n\geq 1.
\end{align}
\end{subequations}
\noindent
We emphasise that equations \eqref{eq:quasiortho-rec}, \eqref{eq:r&l} and \eqref{eq:ortho-conditions} have been adapted from \cite{IsmWan} as we have chosen to consider monic orthogonal and quasi-orthogonal polynomials, whilst the original formulation referred to orthonormal polynomials and their quasi-orthogonal counterparts. 

For a family of quasi-orthogonal polynomials of generic order $(k,r)$ the three-point relation in~\eqref{eq:quasiortho-rec} can be further generalised, and the coefficients of $R_{n+1}, R_n$ and $R_{n-1}$ would be of degree up to $kr,(kr+1)$ and $kr$ in $x$, respectively \cite{chihara} (the degrees being exactly those listed if ${\color{black}\mathcal{c}}_{n,k}$ and ${\color{black}\mathcal{c}}_{n,k-1}$ are non zero). The conditions ensuring the orthogonality of quasi-orthogonal polynomials of degree $(k,r)$ have been recently provided in~\cite{Bracciali}.

\subsection{Skew-orthogonal polynomials}
\label{sec:skewintro}
We now define a class of skew-orthogonal polynomials and summarise some results from~\cite{Adler2000_skew,adler2002}, establishing a map between skew-orthogonal polynomials and suitably defined orthogonal ones. Restricting ourselves to the specific Freud weight $w(x;t)$ in \eqref{eq:freud_weight},   
we introduce the skew-symmetric inner product~$\langle~\cdot~,~\cdot~\rangle_t$ on the real line {\color{black}for suitable functions $f$ and $g$}
\begin{equation}
    \langle\, f,\,g\,\rangle_t = \frac{1}{2}\int_{\mathbb{R}^2} \!f(x)\,g(y)\,  {\rm sgn}(y-x)\,w(x;t)\,w(y;t)\, {\rm d}x\,  {\rm d}y.
    \label{eq:skew_product}
\end{equation}
{\color{black}In particular, we assume $f$ and $g$ to be polynomials.}   
We define the semi-infinite vector of monic polynomials $ Q(x; t) = (\,Q_0(x; t), Q_1(x; t),\, \dots\,)^{\top}$, with elements of the form 
\begin{equation}
    Q_j(x;t) = x^j + b_{j-1}(t)\,x^{j-1}+ \dots + b_{0}(t), \qquad b_j(t) \equiv 1, ~ b_{j-1}(t), \dots, b_0(t) \in \mathbb{R},~\forall~j \in \mathbb{N}_0,
\end{equation}
skew-orthogonal with respect to the product $\langle~\cdot~,~\cdot~\rangle_t$, i.e.\ 
\begin{equation}\label{eq:skew-ortho_monic}
   \langle\, Q_{2i}, Q_{2j+1}\,\rangle_t = -\langle \,Q_{2j+1}, Q_{2i}\,\rangle_t = \delta_{i,j}\,r_j(t) , \qquad 
    \langle\, Q_{2i}, Q_{2j}\,\rangle_t = 0 = \langle\, Q_{2i+1}, Q_{2j+1}\, \rangle_t \, 
\end{equation}
for all $i,\,j\in\mathbb{N}_0$. Moreover, we define $q(x; t) =\left(q_0(x; t), q_1(x; t), \dots \right)^{\top}$ the vector of skew-orthonormal polynomials such that
\begin{equation}
    \langle\, q_{2i}, q_{2j+1}\,\rangle_t = -\langle \,q_{2j+1}, q_{2i}\,\rangle_t = \delta_{i,j}, \qquad 
    \langle\, q_{2i}, q_{2j}\,\rangle_t = 0 = \langle\, q_{2i+1}, q_{2j+1}\, \rangle_t .
\end{equation}
The elements of $Q(x;t)$ and $q(x;t)$ are mapped into each other via the normalisation {\color{black}coefficients} in~\eqref{eq:skew-ortho_monic} $\{r_j(t)\}_{j\in\mathbb{N}}$  such that 
\begin{equation}\label{eq:ortho_monic}
Q_{2j}(x;\, t) = \sqrt{r_j(t)} \,q_{2j}(x;t), \qquad Q_{2j+1}(x;\, t) = \sqrt{r_j(t)} \,q_{2j+1}(x;t), \quad \forall ~ j \in \mathbb{N}_0.
\end{equation}
\begin{remark}
The skew-orthogonal polynomials $Q(x;t)$ in~\eqref{eq:skew-ortho_monic} are not uniquely defined, as their skew-orthogonality structure is not affected by transformations of the type 
\begin{equation}
Q_{2k+1}\rightarrow Q_{2k+1} + a_{2k}\,Q_{2k} ,    \label{eq:Qodd_shift}
\end{equation}
for arbitrary $a_{2k}$ (see e.g.~\cite[eq.\ (2.4)]{Adler2000_skew}). Therefore, in the rest of this work we fix $Q_0(x;t) = 1 ={\color{black}=P_0(x;\,t)}$ and $Q_1(x;t)=x{\color{black}=P_1(x;\,t)}${\color{black}. Moreover, we choose even polynomials $Q_{2n}(x; t)$ as linear combinations of even powers of $x$ only, and $Q_{2n+1}(x; t)$ as a combination of odd powers of $x$ only. We note that it is enough to assume odd parity of polynomials $\{Q_{2j+1}(x;\,t)\}_{j\in\mathbb{N}_0}$ to automatically ensure even parity for $\{Q_{2j}(x;\,t)\}_{j\in\mathbb{N}_0}$.}
\end{remark}
We will see in the following that the skew-symmetric inner product $\langle~\cdot~,~\cdot~\rangle_t$ in~\eqref{eq:skew_product} is deeply related to the symmetric one $(~\cdot~,~\cdot~)_t$ defined as 
\begin{equation}
   \left(\,f,\, g\,\right)_t = \int_{\mathbb{R}}f(x)\,g(x) \,(w(x;t))^2\,  {\rm d}x,
   \label{eq:product}
\end{equation}
{\color{black}i.e.\ with the re-scaled Freud weight~\eqref{eq:omega_freud}}. As per the framework in Section \ref{sec:orthointro}, we introduce the vector $P(x;t)=(P_0(x;t),P_1(x;t),\dots)^{\top}$ of monic polynomials orthogonal with respect to $(~ \cdot~,~\cdot~)_t$, and its orthonormal counterpart $p(x;t)=(p_0(x;t),p_1(x;t),\dots)^{\top}$, such that their elements satisfy
\begin{equation}\label{eq:orthomonic}
   \left(\,P_n, P_m\,\right)_t = \delta_{n,m}\,h_m(t), \qquad  (\,p_n, p_m\,)_t = \delta_{n,m}, \qquad \forall~ n,m \in \mathbb{N}_0.
\end{equation}
Here $\{h_j(t)\}_{j\in\mathbb{N}_0}$ map monic polynomials onto normal ones via $P_i(x;t) = \sqrt{h_i(t)}\, p_i(x; t)$ for all $i\in\mathbb{N}_0$.

The normalisation coefficients $\{h_j(t)\}_{j\in\mathbb{N}_0}$ in~\eqref{eq:orthomonic} and $\{r_j(t)\}_{j\in\mathbb{N}_0}$ in~\eqref{eq:skew-ortho_monic} populate respectively the semi-infinite matrices $\mathcal{D}(t)$ and $\mathcal{R}(t)$, {\color{black} formally defined via the expressions for their elements:} 
\begin{subequations}
\label{eq:definition_R_D_Q}
\begin{align}
\label{eq:diag_entries}
    \mathcal D_{i,j}(t) = \left(\,P_i,\,P_j\,\right)_t &=  \delta_{i,j}\,h_j(t),
    \\[.5ex]
    \mathcal R_{2i,2j+1}(t) =  -\mathcal R_{2j+1, 2i}=\langle\, Q_{2i}, Q_{2j+1}\,\rangle_t &= \delta_{i,j}\,r_j(t),\\[.5ex]
    \mathcal R_{2i,2j}(t) = \mathcal R_{2i+1,2j+1}(t) &=0,
    \end{align}
\end{subequations}
with $i,j\in \mathbb{N}_0$. In particular, $\mathcal{D}(t)$ is a diagonal matrix and $\mathcal{R}(t)$ is constituted by $2\times 2$ skew-symmetric blocks along the main diagonal.

As detailed in a more general framework in \cite{Adler1999ThePL, Adler2000_skew, adler2002}, the connection between orthogonal polynomials $P(x;t)$ and skew-orthogonal polynomials $Q(x;t)$ is realised through a map between the scalar products $\langle~\cdot~,~\cdot~\rangle_t$ and $(~\cdot~,~\cdot~)_t$. In the following, we retrace the steps {\color{black}in~\cite{Adler1999ThePL}} to realise this map, describing the approach discussed there for our {\color{black}specific} choice of weight $w(x;t)$. 

{\color{black} First, we define appropriate semi-infinite matrices and their inverses, as discussed below.} It is an application of {\color{black}\cite[Prop.\ 6.1]{adler2002}}, that for 
{\color{black}all polynomials} $f$ and $g$, {\color{black}the skew-symmetric product can be expressed as the symmetric product involving the operator
$\overline{\mathcal n}^{\,-1}$ as}
\begin{equation}
    (\,f,\, \overline{\mathcal{n}}^{\,-1}\,g\,)_t =  - \langle \,f,\, g \,\rangle_t , 
    \label{eq:2to1bis}
\end{equation}
{\color{black}where $\overline{\mathcal n}^{\,-1}$ denotes the integral operator
\begin{equation}
    \bigl(\overline{\mathcal n}^{\,-1}g\bigr)(x)
    =\frac{1}{2\,w(x;t)}\int_{\mathbb R}\textup{sgn}(x-y)\,w(y;t)\,g(y)\,{\rm d}y.
    \label{eq:int_rep_n}
\end{equation}
This coincides with the inverse of {\color{black} an appropriately defined operator} $\overline{\mathcal n}$ in the sense of~\cite[eq.\ (6.1)]{adler2002},} which 
for the
Freud weight $w(x;t)$ in~\eqref{eq:freud_weight},
reads as\footnote{Namely, for weights of the general form $\exp{(-V(x))}$, one needs $V'(x) = \frac{1}{2}\frac{G(x)}{F(x)}$ with $F(x)$ and $G(x)$ polynomials. In this more general framework, $\overline{\mathcal{n}} = F(x) \,\frac{d}{\text{d}x}+\frac{1}{2}\left({F(x)'-G(x)}\right)$, and the map is between $\langle ~\cdot~,~\cdot~\rangle$ and the symmetric inner product $(~\cdot~,~\cdot~)$ with weight $\exp{(-2V(x) + \log F(x))}$ \cite{Adler1999ThePL,Adler2000_skew, adler2002}.  In our case $F(x) = 1$ and $G(x) = 8 x^3 -4 t x$.} 
\begin{equation}
    \overline{\mathcal{n}} = \frac{ {\rm d}}{ {\rm d}x} - 4 x^3 + 2t x  .
    \label{eq:operator_n}
\end{equation} 
With the operator $\overline{\mathcal{n}}$ we construct the {\color{black} semi-infinite} skew-symmetric matrix $\mathcal{N}(t)$ with elements 
\begin{equation}
   \mathcal{N}_{i,j}(t) = \left(\,P_i, \overline{\mathcal{n}}\, P_j\,\right)_t,\,\qquad i,j\in\mathbb{N}_0.
    \label{eq:matrix_N}
\end{equation}
{\color{black}
The semi-infinite matrices introduced above are defined through their entries and are understood formally: the products are computed entry by entry, and for each fixed pair of indices only finitely many terms contribute. This follows from the structure of the matrices: $\mathcal D$ is diagonal, $\mathcal N$ has four
nonzero diagonals, $\mathcal R$ have two nonzero diagonals each.}

{\color{black} An operator that plays a pivotal role in the following is the semi-infinite transition matrix $\mathcal{Q}(t)$ mapping monic orthogonal polynomials~\eqref{eq:orthomonic} into monic skew-orthogonal polynomials~\eqref{eq:skew-ortho_monic}} 
\begin{equation} 
\label{eq:definition_Q}
Q(x;t) = \mathcal{Q}(t)\,P(x;t). 
\end{equation}
{\color{black}We also }introduce the semi-infinite skew-symmetric matrix $\mathcal M(t)$ such that 
\begin{equation}
    \mathcal M_{i,j}(t) = \left(\,P_i, \overline{\mathcal{n}}^{\,-1} P_j\,\right)_t, \qquad i,j\in\mathbb{N}_0,
\end{equation}
which implies the following identity: 
\begin{equation}
    \mathcal M\, \mathcal D^{-1} \,\mathcal N = \mathcal D.
    \label{eq:decompositionM}
\end{equation} 
This follows straightforwardly from
\begin{equation}\label{eq:prodMDN}
   \left( \mathcal M\, \mathcal D^{-1} \,\mathcal N\right)_{i,j} = -\sum_{k\geq 0}\frac{\left(\,\overline{\mathcal{n}}^{\,-1}P_i,\,  P_k\,\right)_t}{\left(\,P_k,\, P_k\,\right)_t}\left(\, P_k,\,  \overline{\mathcal{n}}\, P_j\,\right)_t = -\left(\,\overline{\mathcal{n}}^{\,-1}\, P_i,\, \overline{\mathcal{n}}\,P_j\,\right)_t = \mathcal{D}_{i,j},
\end{equation}
exploiting the decomposition 
\begin{equation}
\overline{\mathcal{n}}\,P_j = \sum_{k\geq 0}\frac{\left(\,P_k,\, \overline{\mathcal{n}}P_j\,\right)_t}{\left(\,P_k,\, P_k\, \right)} P_k.
\label{eq:n_decomposition}
\end{equation}
{\color{black} We emphasise that the product~\eqref{eq:prodMDN} is well defined, {\color{black} as for each pair of indices $i$ and $j$ only finitely many terms appear in the sum}. Although $\mathcal{M}$ is non sparse with entries given by $\mathcal M_{i,j}=-\langle P_i,P_j\rangle_t$ via the mapping~\eqref{eq:2to1bis}, 
{\color{black} on the other hand} we have $\mathcal D^{-1}=\text{diag}\bigl(\left\{h_j(t)^{-1}\right\}_{j \in \mathbb N_0}\bigr)$ with $h_j(t)=(P_j,P_j)_t>0$, and $\mathcal{N}$ has four nonzero diagonals.}

Moreover, it follows from the definition of $\mathcal R(t)$ and $\mathcal Q(t)$ that the {\color{black}following decomposition holds}
\begin{equation}
    \mathcal R = -\mathcal Q \,\mathcal M\, \mathcal Q^\top.
    \label{eq:decompositionR}
\end{equation}
Therefore from equations~\eqref{eq:decompositionM} and \eqref{eq:decompositionR} with some algebra it follows that
\begin{equation}
    \mathcal Q^{\top}\, \mathcal R^{-1} \,\mathcal Q = -\mathcal D^{-1}\,\mathcal N\, \mathcal D^{-1}.
    \label{eq:relation_for_Q}
\end{equation}
{\color{black}The matrix $\mathcal{R}^{-1}$ is well defined, being a diagonal block matrix with $2 \times 2$ blocks of the form 
\begin{equation*}
    \begin{pmatrix}
        0 & -\frac{1}{r_j(t)} \\
        \frac{1}{r_j(t)} & 0
    \end{pmatrix} ,
\end{equation*}
corresponding to $-{r_j(t)}^{-1} \,J$, where $J$ is the $2\times 2$ skew-identity ($J^2=-I$) and $r_j \neq 0$ are the normalisation coefficients in~\eqref{eq:ortho_monic}.}
{\color{black}
This proves the following \cite{Adler2000_skew}\footnote{\color{black} In \cite{Adler2000_skew}, the statement of \ref{prop:decomposition} has been shown for a broader class of weights which include~\eqref{eq:freud_weight}.  }
\begin{prop}\label{prop:decomposition}
   The semi-infinite matrix $\mathcal{N}(t)$ admits the following decomposition
\begin{equation}
    \mathcal{N} = - \mathcal D\,\mathcal{Q}^{\top}\, \mathcal R^{-1}\, \mathcal{Q}\,\mathcal D = -(\mathcal{Q}\,\mathcal D)^{\top}\, \mathcal R^{-1}\, \mathcal{Q}\,\mathcal D.
    \label{eq:Ndecomposition}
\end{equation}
Here $\mathcal{D}(t)$ and $\mathcal{R}(t)$ are the matrices defined in~\eqref{eq:definition_R_D_Q} and $\mathcal{Q}(t)$ is defined in \eqref{eq:definition_Q}.
\end{prop}
}

Note that explicit expressions for the elements of $\mathcal{Q}(t)$ would allow us to generate the skew-orthogonal polynomials \eqref{eq:skew-ortho_monic} in terms of orthogonal ones \eqref{eq:orthomonic}. This is going to be our goal for Section \ref{sec:transition}. To accomplish it, we are going to rely on the properties of orthogonal polynomials \eqref{eq:orthomonic} detailed in Section~\ref{sec:orthointro}.

\begin{remark}\label{remark:Hermitian&Symmetric}
The normalisation coefficients $\{h_i(t)\}_{i\in\mathbb{N}_0}$ and $\{r_i(t)\}_{i\in\mathbb{N}_0}$ can be written in terms of partition functions of certain random matrix ensembles. 
Indeed, following \cite{vanMoerbekenotes} and references therein, $\{h_i(t)\}_{i\in\mathbb{N}_0}$ admit a representation in terms of Hankel determinants of the moment matrix built for powers of $x^k$ via the symmetric product $(~\cdot~,~\cdot~)_t$:
\begin{equation}
    h_n(t) = \frac{\tau^{\textup{Herm}}_{n+1}(t)}{\tau^{\textup{Herm}}_n(t)}\, \quad \text{ with } 
    \tau^{\textup{Herm}}_n(t) = \det( m_n(t)), \quad (m_n(t))_{i,j}=(\,x^i, x^j\,)_t,\;\; 0\leq i,j<n.
\end{equation}
In the context of the random matrix model, $\tau^{\text{Herm}}_n(t)$ corresponds to the partition function for the ensemble of random $n\times n$ Hermitian matrices with weight $w(x;t)^2$, which can be explicitly expressed as 
\begin{equation}
    \tau_n^{\textup{Herm}}(t) = \frac{1}{n!}\int_{\mathbb{R}^n} \prod_{i< j} (z_i-z_j)^2 \, \prod_{i=1}^n w(z_i;t)^2\, {\rm d}z_i. 
\end{equation}
Here the integration variables $\{z_i\}_{i=1}^n$ represent the eigenvalues of the ensemble matrices.
Analogously, the coefficients $\{r_j(t)\}_{j \in \mathbb{N}_0}$ can be formulated in terms of the pfaffian of  the moment matrix $m_n(t)$ built using the skew-symmetric product $\langle\,\cdot,\,\cdot\,\rangle_t$, i.e.\ 
\begin{equation}
    r_n(t) = \frac{\tau^{\textup{Symm}}_{2(n+1)}(t)}{2 \tau^{\textup{Symm}}_{2n}(t)}\, \quad \text{ with } 
    \tau^{\textup{Symm}}_{{\color{black}2}n}(t) = 2^n\,\text{pf}( m_n(t)), \quad (m_n(t))_{i,j}=\langle\,x^i, y^j\,\rangle_t,\;\; 0\leq i,j<n,
    \label{eq:r&tau}
\end{equation}
where the pfaffian for matrices of even dimension is defined as the square root of the determinant. 
Furthermore, the function $\tau_{2n}^{\textup{Symm}}(t)$ constitutes the partition function for the ensemble of random symmetric $2n \times 2n$ matrices with weight $w(x;t)$
\begin{equation}
    \tau_{2n}^{\textup{Symm}}(t)= \frac{1}{(2n)!}\int_{\mathbb{R}^{2n}}\prod_{i<j} |z_i-z_j| \, \prod_{i=1}^{\color{black}2n} w(z_i;t) \, {\rm d}z_i,
    \label{eq:tau_symm}
\end{equation}
with $\{z_i\}_{i=1}^{\color{black}2n}$ representing the eigenvalues of the matrices of the ensemble. {\color{black} We remark that the relationship between skew-orthogonal polynomials and the symmetric matrix ensemble has been established for even-sized matrices only, as per equations \eqref{eq:r&tau} and \eqref{eq:tau_symm}.}
{\color{black} The even-sized} symmetric matrix ensemble has been recently investigated in \cite{BDM1,BDM2} exploiting its relationship with integrable systems and skew-orthogonal polynomials.
\end{remark}

\section{Skew-orthogonal polynomials as quasi-orthogonal ones}\label{sec:transition}
In this section we focus on the mapping between Freud orthogonal polynomials $P(x;t)$ in~\eqref{eq:orthomonic} and their skew-symmetric counterparts $Q(x;t)$ in~\eqref{eq:skew-ortho_monic}.
Our goal is to build upon \cite{Adler2000_skew} and write explicitly the matrix elements of $\mathcal{N}(t)$ in \eqref{eq:matrix_N} exploiting \eqref{eq:3pts} and \eqref{eq:derivative}, in order to find the transition matrix~$\mathcal{Q}(t)$ via the relation \eqref{eq:relation_for_Q}. This mapping will reveal the underlying structure of skew-orthogonal polynomials as appropriate quasi-orthogonal ones.

\subsection{Transition matrix for Freud weights}
The transition matrix $\mathcal{Q}(t)$
mapping {\color{black}orthogonal polynomials $P(x;t)$ into skew-orthogonal ones $Q(x;t)$} satisfies the identity~\eqref{eq:relation_for_Q}, 
with $\mathcal{R}(t)$ and $\mathcal{D}(t)$ defined in equation \eqref{eq:definition_R_D_Q}, and $\mathcal{N}(t)$ introduced in~\eqref{eq:matrix_N}.
As mentioned above, our goal is to explicitly evaluate the matrix elements of $\mathcal{Q}(t)$ exploiting this relationship, which we reformulate as follows for ease of calculation \begin{equation}
\mathcal{R}^{-1}\mathcal{Q} = -(\mathcal{Q}^\top)^{-1}\mathcal{D^{-1}}\,\mathcal{N}\,\mathcal{D^{-1}}.
    \label{eq:relation_for_Q2}
\end{equation}

To start, we focus on the matrix $\mathcal{N}(t)$ and its elements. From the definition of the operator $\overline{\mathcal{n}}$ in~\eqref{eq:operator_n}, it is clear that to formulate $\mathcal{N}(t)$ explicitly we need to express $x^3 P_n(x;t)$ for a generic $n\in\mathbb{N}_0$ as a linear combination of polynomials $P(x;t)$. Therefore we apply~\eqref{eq:3pts} three times and find
\begin{equation}
    x^3\, P_n = P_{n+3} + \left(\beta_{n+2} + \beta_{n+1} + \beta_n\right)P_{n+1} +\beta_n\left(\beta_{n+1}+\beta_n + \beta_{n-1}\right)  P_{n-1} +  \beta_n\beta_{n-1}\beta_{n-2}\,P_{n-3},
\end{equation}
with the understanding that, for $n<3$, quantities labelled with negative indices are null. 
Thus, the action of the operator $\overline{\mathcal{n}}$ is given by 
\begin{equation}
\begin{split}
    \overline{\mathcal{n}}\, P_n &= 4\, \beta_n \beta_{n-1}\beta_{n-2}\, P_{n-3} + \big(n  +2\beta_n\big(t-2(\beta_{n+1}+\beta_n+\beta_{n-1})\!\Big)\big)P_{n-1}\\
    &~~+2\big(t-2(\beta_n + \beta_{n+1}+\beta_{n+2})\!\Big)P_{n+1} - 4\, P_{n+3}, 
    \end{split}
   \label{eq:nPn}
\end{equation}
where we have used the recurrence relation for $\beta_n(t)$ in~\eqref{eq:recurrence_beta}. 
Equation~\eqref{eq:nPn} can be formulated more succinctly using the definition~\eqref{eq:beta_def} and again the relation~\eqref{eq:recurrence_beta}:
\begin{equation}
     \overline{\mathcal{n}}\, P_n = 4\, \frac{h_n}{h_{n-3}}\, P_{n-3} + \frac{n}{2}\,P_{n-1} -\frac{n+1}{2}\,\frac{h_n}{h_{n+1}}\,P_{n+1}  - 4\, P_{n+3}.
    \label{eq:nPn_simple}
\end{equation}
Recalling the decomposition of $\bar{\mathcal{n}}\, P_n$ according to \eqref{eq:n_decomposition}, i.e.
$$
\overline{\mathcal{n}}\,P_n = \sum_{\ell \geq 0}\frac{\left(\,P_\ell\, , \overline{\mathcal{n}}\, P_n\,\right)_t}{\left(\,P_\ell, P_\ell \,\right)_t}\,P_\ell  , 
$$
we have from \eqref{eq:nPn_simple} that $\left(\,P_\ell, \overline{\mathcal{n}}\, P_n\,\right)$ is non-zero only for $\ell \in \{n-3, n-1, n+1, n+3\}$ due to the orthogonality of the polynomials. Therefore the elements of the matrix $\mathcal{N}(t)$ are
\begin{equation}\label{eq:matrix_N_explicit}
\mathcal{N}_{i,j} = \left(P_i, \overline{\mathcal{n}}\, P_j\right)_t =4 \,h_{i+3}\,\delta_{j,i+3}+\frac{i+1}{2}\, h_{i}\, \delta_{j,i+1}-\frac{i}{2} \,h_{i-1}\,\delta_{j,i-1}-4\, h_{i}\, \delta_{j,i-3}.
\end{equation}
 Notice that the only non-zero elements {\color{black}$\mathcal{N}_{i,j}$ are those corresponding to the indices $(i,j)=(2k+1,2k)$, $(i,j)=(2k-1,2k)$ and $(i,j)=(k,k+3)$ }
\begin{equation}
\label{eq:Nelements}
        \mathcal{\mathcal{N}}_{2k+1,2k}=-\frac{(2k+1)}{2}\,h_{2k},\quad 
       \mathcal{N}_{2k-1, 2k}= k\, h_{2k-1} , \quad
       \mathcal{N}_{k,k+3}=4 \,h_{k+3} ,
\end{equation}
and their counterparts with indices in the opposite order, as $\mathcal{N}_{i,j}=-\mathcal{N}_{j,i}$. We thus have an explicit expression for the right hand side of~\eqref{eq:relation_for_Q}  for a generic matrix element in terms of the coefficients $\{h_j(t)\}_{j\in\mathbb{N}_0}$ only:
\begin{equation}
    \left(-\mathcal{D}^{-1}\,\mathcal{N}\,\mathcal{D}^{-1}\right)_{i,j}=-\frac{4}{h_{i}}\, \delta_{j,i+3}-\frac{i+1}{2 h_{i+1}}\,  \delta_{j,i+1}+\frac{i}{2h_i}\, \delta_{j,i-1}+\frac{4}{h_{j-3}}\, \delta_{j,i-3}, 
\end{equation}
{\color{black}with $h_j \neq 0$ for all $j \ge 0$, being normalisation coefficients for orthogonal polynomials. }

We now consider the transition matrix $\mathcal{Q}(t)$. In particular, we can characterise its elements as 
\begin{equation}
    \mathcal{Q}_{i,j}(t) = \delta_{i,j} + \vartheta(i-j)\, \mathcal{q}_{i,j}(t), 
\end{equation}
with $\vartheta(x)$ being the Heaviside function giving $1$ for $x>0$ and $0$ otherwise. The question is then evaluating the functions~$\mathcal{q}_{i,j}(t)$. 
We introduce the functions $\varphi,\sigma \colon \mathbb{N}_0 \to \mathbb{Z}_2$ as in~\cite{BDM2}
\begin{equation}
    \varphi(i) = \begin{cases}
        1 \quad i \text{ even} \\
        0 \quad i \text{ odd}
    \end{cases}, \qquad  \sigma(i) = \begin{cases}
        1 \quad i \text{ odd} \\
        0 \quad i \text{ even}
    \end{cases}, 
\end{equation}
with the following properties in terms of their arguments: 
    \begin{align*}
    \varphi(i+n)&=\begin{cases}
        \varphi(i), \quad |n| \text{ even} \\
        \sigma(i),\, \quad |n| \text{ odd}
    \end{cases}, \qquad \sigma (i+n)=\begin{cases}
        \sigma(i), \quad |n| \text{ even}  \\
        \varphi(i), \quad |n| \text{ odd}
    \end{cases} , \\
    \begin{split}   &\varphi(i+j)=\varphi(i)\,\varphi(j)+\sigma(i)\,\sigma(j), \qquad \sigma(i+j)= \varphi(i)\,\sigma(j)+\sigma(i)\,\varphi(j) , \\[1.5ex]
    &\varphi(i)^{m}=\varphi(i), \qquad \sigma(i)^{m}=\sigma(i),  \qquad \varphi(i)\,\sigma(i)=0, \quad m \in \mathbb{N}.
\end{split} 
\end{align*}
Hence, $\mathcal{Q}(t)$ can be further characterised as 
\begin{equation}\label{eq:Qmatrix_structure}
\mathcal{Q}_{i,j} = \delta_{i,j} + \vartheta(i-j)\big(\sigma(i)\,\sigma(j) + \varphi(i)\,\varphi(j)\!\Big)\,\mathcal{q}_{i,j},
\end{equation}
meaning that the non-zero values of $\mathcal{q}_{i,j}(t)$ are those for which either both the indices are even numbers or they are both odd. 

As the inverse of the coefficient matrix $\mathcal{R}(t)$ appears in the left hand side of~\eqref{eq:relation_for_Q2} as well as $\mathcal{Q}(t)$, it is helpful to provide the generic element of $\mathcal{R}^{-1}(t)$, {\color{black} which reads as} 
\begin{equation}\label{eq:Rinverse_element}
    (\mathcal{R}^{-1})_{i,j} = \sigma(i)\,\varphi(j)\, \delta_{i, j+1}\, \frac{1}{r_{\frac{j}{2}}}-\varphi(i)\,\sigma(j)\, \delta_{j,i+1}\,\frac{1}{r_{\frac{i}{2}}} ,
\end{equation}
{\color{black}with $r_j \neq 0$ for all $j \ge 0$, being normalisation coefficients for skew-orthogonal polynomials.}
Therefore {\color{black} each element of }the left hand side of \eqref{eq:relation_for_Q2} {\color{black}is given by the sum of finitely many terms, as $\mathcal Q$ in \eqref{eq:Qmatrix_structure} has only
finitely many nonzero elements in each row and column, and} becomes
\begin{equation}
    \begin{split}
\left(\mathcal{R}^{-1}\,\mathcal{Q}\right)_{i,j} &= \sigma(i)\varphi(j)\left(\frac{\delta_{j,i-1} }{r_{\frac{i-1}{2}}}+\vartheta\left(i-(j+1)\right)\frac{\mathcal{q}_{i-1,j}}{r_{\frac{i-1}{2}}}\right) - \varphi(i)\sigma(j)\left(\frac{\delta_{j,i+1} }{r_{\frac{i}{2}}}+\vartheta\left(i-(j-1)\right)\frac{\mathcal{q}_{i+1,j}}{r_{\frac{i}{2}}}\right).
    \label{eq:Rinv_Q}
    \end{split}
\end{equation}

We can now turn to {\color{black}the right hand side of} \eqref{eq:relation_for_Q2}.
Due to the structure of $\mathcal{Q}(t)$, $\mathcal{D}^{-1}(t)$, and $\mathcal{N}(t)$ {\color{black}as per equations~\eqref{eq:q_ij_explicit}, \eqref{eq:diag_entries} and \eqref{eq:matrix_N_explicit}} respectively, it can be easily shown that 
\begin{equation}
    \left(-(\mathcal{Q}^\top)^{-1}\mathcal{D}^{-1}\mathcal{N}\mathcal{D}^{-1}\right)_{i,j} = 0, \qquad \text{  for } i\geq j+4,
\end{equation}
i.e.\ all the elements in diagonals which are lower than the forth below the principal diagonal are null. {\color{black}Here $(\mathcal Q^\top)^{-1}$ is well defined and explicit: since
$\mathcal Q^\top$ is unit upper triangular, its inverse is determined recursively from the top-left corner, and each entry is given by a finite sum. Moreover, the strictly upper triangular part of $\mathcal Q^\top$ has non vanishing entries only for $j-i=2$ in even rows and $j-i=2$ or $4$ in odd rows (see \eqref{eq:q_null} and \eqref{eq:q_ij_explicit}). Therefore, each row and column contains at most two non vanishing entries, and then every product and sum appearing in the expression above is finite.}
From identity~\eqref{eq:relation_for_Q2} and expression \eqref{eq:Rinv_Q} it follows that
\begin{subequations}
\begin{align}
    \mathcal{q}_{2(\ell+k), 2\ell} &= 0, \qquad \text{ if } k\geq2,\\
    \mathcal{q}_{2(\ell+k)+1, 2\ell+1} &= 0, \qquad \text{ if }  k\geq3,
    \label{eq:q_null}
\end{align}
\end{subequations}
and therefore the expression in~\eqref{eq:Rinv_Q}, i.e.\ the left hand side of~\eqref{eq:relation_for_Q2},  becomes: 
\begin{equation}
\begin{split}
    \hspace*{-1ex}\left(\mathcal{R}^{-1}\mathcal{Q}\right)_{i,j} &=\sigma(i)\!\left(\frac{\delta_{j,i-1}}{r_\frac{i-1}{2}}+\frac{\mathcal{q}_{i-1,i-3}}{r_{\frac{i-1}{2}}}\delta_{j,i-3}\right)\!- \varphi(i)\!\left(\frac{\delta_{j,i+1}}{r_{\frac{i}{2}}} +\frac{\mathcal{q}_{i+1,i-1}}{r_{\frac{i}{2}}}\delta_{j,i-1}+\frac{\mathcal{q}_{i+1,i-3}}{r_{\frac{i}{2}}}\delta_{j,i-3}\right)\!.   
    \label{eq:stepR}
\end{split} 
\end{equation}
The right hand side of~\eqref{eq:relation_for_Q2} is instead
{\color{black} 
\begin{equation}
    \begin{split}
        \left(-(\mathcal{Q}^\top)^{-1}\mathcal{D}^{-1}\mathcal{N}\mathcal{D}^{-1}\right)_{i,j} &=  \left( \frac{4}{h_{i+1}}\mathcal{q}_{i+2,i}\mathcal{q}_{i+4,i+2} - \frac{i+2}{h_{i+2}}\mathcal{q}_{i+2,i} - \frac{4 \sigma(i)}{h_{i+1}}\mathcal{q}_{i+4,i}
        -\frac{i+1}{2h_{i+1}}  \right)\delta_{j,i+1} \\
        &\quad  - \left(\frac{4 \mathcal{q}_{i+2,i}}{h_{i-1}} - \frac{i}{2 h_{i}}\right)\delta_{j,i-1}+\frac{4}{h_{i-3}}\delta_{j,i-3} +\mathcal{A}_{i,j},
   \label{eq:stepD}
    \end{split}
\end{equation}}

\noindent
where $\mathcal{A}_{i,j}$ is in non-zero only for $j = i + (2m+1)$, $m \in \mathbb{N}$. This {\color{black}element} will not play a role in our analysis, and  we will remark briefly on it in the next section. We have thus characterised the structure of the matrix equation~\eqref{eq:relation_for_Q2}. In the next section we exploit~\eqref{eq:stepR} and~\eqref{eq:stepD} to determine $\mathcal{q}_{i,j}(t)$ explicitly in terms of $\{h_j(t)\}_{j\in\mathbb{N}_0}$ and $\{r_{j}(t)\}_{j\in\mathbb{N}_0}$.

\subsection{The map between orthogonal and skew-orthogonal polynomials}
Recall that, for \eqref{eq:relation_for_Q2} to be satisfied, expressions \eqref{eq:stepR} and \eqref{eq:stepD} need to coincide. We obtain the following identities for the elements $\mathcal{q}_{i,j}(t)$
\begin{subequations}\label{eq:q_ij_explicit}
\begin{align}
\label{eq:xi}
    \xi_\ell(t) &:= \mathcal{q}_{2\ell+2, 2\ell}(t)= 4\,\frac{ r_{\ell+1}(t)}{h_{2\ell}(t)} , \\
    \label{eq:psi} 
   \psi_\ell(t) &:= \mathcal{q}_{2\ell+3,2\ell+1}(t) = 4\,
   \xi_{\ell+1}(t)\,   \frac{r_{\ell+1}(t)}{h_{2\ell+1}(t)}-(\ell+1)\,\frac{r_{\ell+1}(t)}{h_{2\ell+2}(t)} ,\\
\label{eq:zeta}
\zeta_\ell(t) &:= \mathcal{q}_{2\ell+5, 2\ell+1}(t)= -4\, \frac{r_{\ell+2}(t)}{h_{2\ell+1}(t)} ,
\end{align}
\end{subequations}
where we introduce the coefficients $\xi_\ell(t)$, $\psi_\ell(t)$ and $\zeta_\ell(t)$ with $\ell\in\mathbb{N}_0$, to enhance the legibility in the following. {\color{black} All the elements are well defined, since $h_j \neq 0$ for all $j \ge 0$.}
We remark that \eqref{eq:xi} follows from the comparison of matrix elements $(2\ell +3, 2\ell)$ in identity~\eqref{eq:relation_for_Q2}, whilst \eqref{eq:psi} follows from elements $(2\ell+2, 2\ell+1)$ and \eqref{eq:zeta} from elements $(2\ell+4, 2\ell+1)$.
Moreover, from \eqref{eq:Qmatrix_structure} and \eqref{eq:q_null} recall that the only non-zero elements of $\mathcal{Q}(t)$ are given by \eqref{eq:q_ij_explicit} and diagonal ones $\mathcal{Q}_{i,i} = 1$.

Based on the definition of the coefficients $\{\beta_j(t)\}_{j \in \mathbb{N}_0}$ as per equation~\eqref{eq:beta_def}, we note that \eqref{eq:psi} and \eqref{eq:zeta} can be expressed in terms of $\{\beta_j(t)\}_{j \in \mathbb{N}_0}$ and $\{\xi_j(t)\}_{j \in \mathbb{N}_0}$ as
\begin{subequations}
\begin{align}
   \psi_\ell &= \frac{\xi_\ell\,\xi_{\ell+1}}{\beta_{2\ell+1}}-\frac{(\ell+1)\,{\color{black}\xi_\ell}}{4\,\beta_{2\ell+1}\,\beta_{2\ell+2}},\label{eq:recurrence_psi}\\
  \zeta_\ell & = - \beta_{2\ell+2}\,\xi_{\ell+1}.\label{eq:recurrence_zeta}
\end{align}
 \label{eq:recurrence_psi&zeta}
\end{subequations}
{\color{black}
The coefficients $\psi_\ell$ are well defined for $\ell \ge 0$, as its denominator {\color{black}features non vanshing terms} $\beta_j{\color{black}(t)}$ with $j>0$.}
Moreover, notice that, at this stage, we do not have a direct way to calculate the coefficients $\{r_j(t)\}_{j\in\mathbb{N}_0}$, whilst the $\{h_j(t)\}_{j \in \mathbb{N}_0}$ are  given in terms of the recurrence coefficients $\{\beta_j(t)\}_{j \in \mathbb{N}_0}$ as per equation~\eqref{eq:beta_def}. This issue can be tackled by considering one more relation following from matrix elements $(2\ell, 2\ell+1)$ in the identity~\eqref{eq:relation_for_Q2}, i.e.
\begin{equation}
    1 = \frac{(2\ell+1)}{2}\, \frac{r_\ell}{h_{2\ell+1}}+ (\ell+1)\,\frac{r_\ell}{h_{2\ell+2}}\, \mathcal{q}_{2(\ell+1),2\ell}-4\, \frac{r_\ell}{h_{2\ell+1}}\, \mathcal{q}_{2(\ell+2),2(\ell+1)}\,\mathcal{q}_{2(\ell+1){\color{black},}\,2\ell}.
    \label{eq:recxi_withr}
\end{equation}
{\color{black}We now exploit some previous equations and use some algebra to express this in a more compact form. First note that, according to \eqref{eq:xi}, $\mathcal{q}_{2\ell +2, 2\ell}(t)$ corresponds to $\xi_\ell(t)$ and we can reformulate $\{r_j(t)\}_{j\in\mathbb{N}_0}$ in terms of $\{\xi_j(t)\}_{j\in\mathbb{N}_0}$, i.e. $r_\ell(t) = h_{2\ell-2}\xi_{\ell-1}(t)/4$. This gives
\begin{equation*}
    1 = \frac{(2\ell+1)}{2}\, \frac{h_{2\ell-2}\xi_{\ell-1}}{4h_{2\ell+1}}+ (\ell+1)\,\frac{h_{2\ell-2}\xi_{\ell-1}}{4 h_{2\ell+2}}\, \xi_\ell-\, \frac{h_{2\ell-2}\xi_{\ell-1}}{h_{2\ell+1}}\, \xi_{\ell+1}\,\xi_\ell.
\end{equation*}
Recalling \eqref{eq:beta_def} in order to express the coefficients $\{h_{i}(t)\}_{i\in\mathbb{N}_0}$ in terms of $\{\beta_{i}(t)\}_{\mathbb{N}_0}$ and rearranging the formula above we find}
the following recurrence equation for $\{\xi_j(t)\}_{j\in\mathbb{N}_0}$
\begin{equation}
 4 \,\beta_{2\ell-1} \,\beta_{2\ell}\,\beta_{2\ell+1}-\frac{2\ell+1}{2}\,\xi_{\ell-1} - (\ell + 1)\,\frac{\xi_{\ell-1} \, \xi_{\ell}}{\beta_{2\ell+2}\,} + 4\,\xi_{\ell-1} \,\xi_\ell \,\xi_{\ell+1}   =0, \quad \ell \ge {\color{black}1}.
 \label{eq:recurrence_xi}
\end{equation}
{\color{black} Explicit expressions for initial conditions $\xi_0(t)$ and $\xi_1(t)$ are provided in Appendix \ref{app:initial}.} 
We speculate that this might constitute a version of a discrete Painlev\'e (see \cite{Sakai2001,GraRam2004}) in higher dimensions. This will be the object of further investigation in future works. 

Recall that the parameters $\{\beta_j(t)\}_{j\in\mathbb{N}_0}$ satisfy the closed three-point relation \eqref{eq:recurrence_beta}. Therefore,~\eqref{eq:recurrence_xi} can be coupled with \eqref{eq:recurrence_beta} to find $\{\xi_j(t)\}_{j\in\mathbb{N}_0}$. This in turn allows us to find $\{\psi_j(t)\}_{j\in\mathbb{N}_0}$ and $\{\zeta_j(t)\}_{j\in\mathbb{N}_0}$ according to \eqref{eq:recurrence_psi} and \eqref{eq:recurrence_zeta}. As a side note, we observe by direct inspection that the elements $\mathcal{A}_{i,j}$ introduced in \eqref{eq:stepR} which would be in principle non-zero are null modulo relations \eqref{eq:q_ij_explicit}, consistently with identity \eqref{eq:relation_for_Q2} and expressions \eqref{eq:stepR} and \eqref{eq:stepD}.

To summarise, the parameters $\xi_\ell(t)$, $\zeta_\ell(t)$, $\psi_\ell(t)$ with $\ell \in \mathbb{N}_0$ completely characterise the transition matrix $\mathcal{Q}(t)$ in~\eqref{eq:definition_Q}. The recurrence relations for $\beta_\ell(t)$ and $\xi_\ell(t)$ in~\eqref{eq:recurrence_beta} and~\eqref{eq:recurrence_xi} provide us with explicit formulas we can use to evaluate the coefficients for any given value of $t$. 
We have thus proved the following theorem.
\begin{theorem}\label{thm:skew_quasi}
The monic polynomials $Q(x;t) =(Q_0(x;t),Q_1(x;t),\dots)^\top$ skew-orthogonal with respect to the inner product \eqref{eq:skew_product} can be expressed as linear combinations of monic orthogonal polynomials with respect to \eqref{eq:product} as
\begin{subequations}\label{eq:QandP}
\begin{align}
        Q_{2n}(x; t) &= P_{2n}(x; t) + \xi_{n-1}(t) \,P_{2n-2}(x; t),\qquad \forall~ n \in\mathbb{N}_{\color{black}0}, \label{eq:evenQandP}\\
        Q_{2n+1}(x; t) &= P_{2n+1}(x; t) + \psi_{n-1}(t) \,P_{2n-1}(x; t) + \zeta_{n-2}(t)\,P_{2n-3}(x; t),\qquad \forall ~ n \in\mathbb{N}_{\color{black}0}. \label{eq:oddQandP}
\end{align}
\end{subequations}
The coefficients $\xi_n(t)$, $\psi_n(t)$ and $\zeta_{n}(t)$ satisfy the relationships \eqref{eq:recurrence_xi}, \eqref{eq:recurrence_psi} and \eqref{eq:recurrence_zeta} respectively. In the expressions above, all quantities are null when indexed with negative indices. 
\end{theorem}
The form of the first elements of the skew-orthogonal polynomials $Q(x;t)$ with even and odd degrees in terms of the coefficients $\{\beta_{i}(t)\}_{i \in \mathbb{N}_0}$ and $\{\xi_{i}(t)\}_{i \in \mathbb{N}_0}$ are listed in Appendix~\ref{app:formskew}. 

\begin{figure}[ht]
\centering
    \includegraphics[width=1\textwidth]{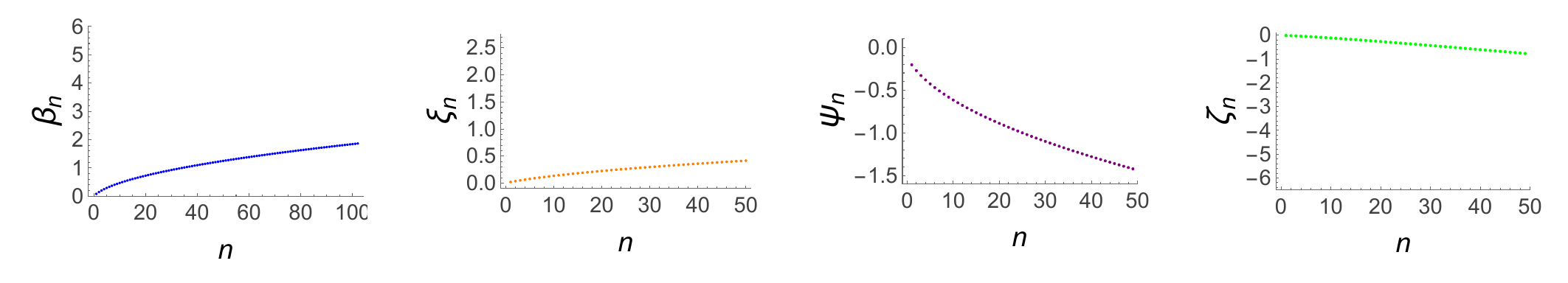} \\[-2.5ex]
    $t=-2.5$ \\[1.5ex]
    \includegraphics[width=1\textwidth]{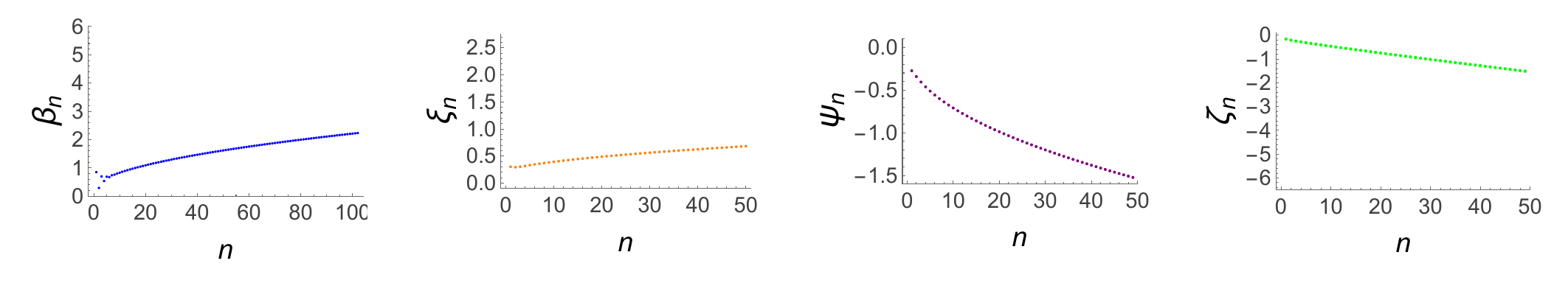} \\[-2.5ex]
    $t=2$ \\[1.5ex]
    \includegraphics[width=1\textwidth]{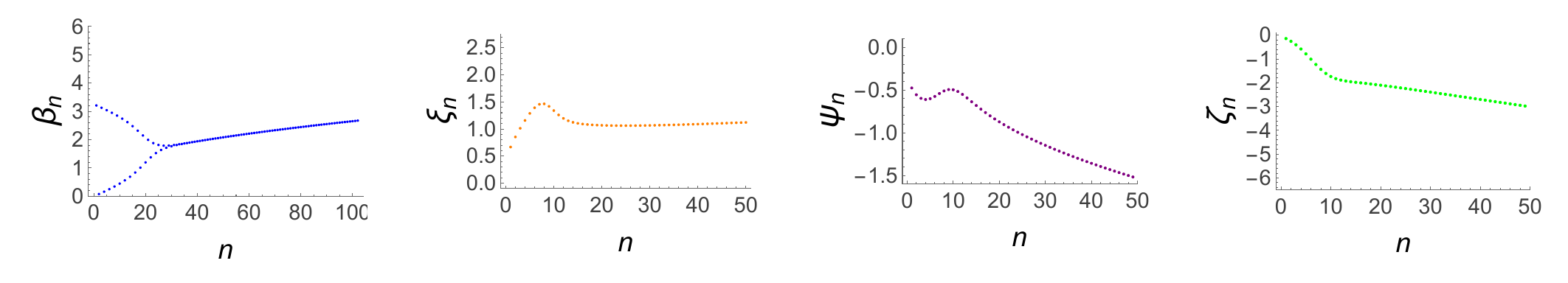} \\[-2.5ex]
    $t=6.5$ \\[1.5ex]
    \includegraphics[width=1\textwidth]{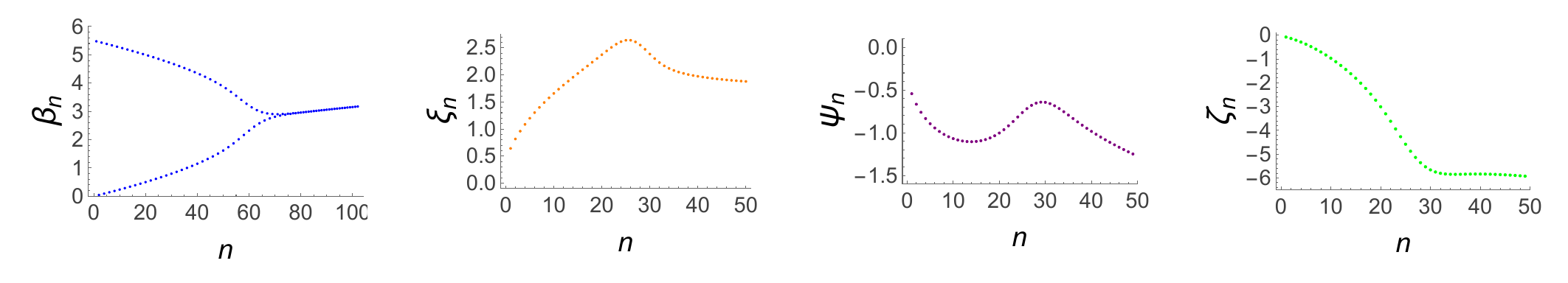} \\[-2.5ex]
    $t=11$ 
\caption{From left to right: $\beta_n(t)$ {\color{black} (blue dots)}, $\xi_n(t)$ {\color{black} (blue dots)}, $\zeta_n(t)$ {\color{black} (purple dots)} and $\psi_n(t)$  {\color{black} (green dots)} shown as functions of $n$ for increasing values of the parameter $t$ in the Freud weight \eqref{eq:freud_weight}.}
\label{fig:plots}
\end{figure}

We note, incidentally, that identities \eqref{eq:recurrence_psi&zeta} and  \eqref{eq:recurrence_xi} can be in principle solved recursively for any fixed value of the parameter $t$. This first requires tackling the recursive relation \eqref{eq:recurrence_beta} in order to have suitable values for the parameters $\{\beta_i(t)\}_{i\in\mathbb{N}_0}$ which can then be used within  \eqref{eq:recurrence_xi}, \eqref{eq:recurrence_psi} and \eqref{eq:recurrence_zeta}. Indeed, we emphasise that it is a key property of equations \eqref{eq:QandP} that all the parameters involved 
can be -- more or less directly -- expressed in terms of the variables $\{\beta_{i}(t)\}_{i\in\mathbb{N}_0}$, as they appear as coefficients of the recursive relation determining $\{\xi_i(t)\}_{i\in\mathbb{N}_0}$ \eqref{eq:recurrence_xi}. Therefore, the mapping between skew-orthogonal and orthogonal polynomials is underpinned by the variables $\{\beta_{i}(t)\}_{i\in\mathbb{N}_0}$ and their recurrence relations~\eqref{eq:recurrence_beta}. However, it should be remarked that the starting values $\xi_0(t)$ and $\xi_1(t)$ need {\color{black} nonetheless} to be evaluated separately, as detailed in Appendix \ref{app:initial}.
We include plots for $\beta_n(t)$, $\xi_n(t)$, $\zeta_n(t)$ and $\psi_n(t)$ as functions of $n$ for a selection of values of $t$ in Figure~\ref{fig:plots}. {\color{black}  We point out that similar plots for $\{\beta_i(t)\}_{i\in\mathbb{N}_0}$ in a more general framework have been provided elsewhere in the literature, for example most recently in \cite{HermitianPRE, Clarkson25}. The two curves appearing for $\beta_n(t)$ as function of $n$ given high enough values of $t$ correspond to even and odd values of $n$ respectively, which then converge to the same curve for $n$ large enough. This behaviour has been investigated in e.g. \cite{HermitianPRE, Clarkson25}. Similar investigation of $\xi_n(t)$, $\psi_n(t)$ and $\zeta_n(t)$ as functions of both $n$ and $t$ will be object of future work.
} 

We note that equations \eqref{eq:QandP} show how monic skew-orthogonal polynomials resemble quasi-orthogonal polynomials with respect to the the scalar product \eqref{eq:product} -- {\color{black}the even-degree ones are quasi-orthogonal} of order $(1,2)$ and {\color{black}the odd-degree ones are quasi-orthogonal} of order $(2,2)$. 
Due to this structure, it follows from \eqref{eq:quasiortho_int} that
\begin{equation}
\begin{split}
   \int Q_{2n}(x;t)\,Q_{2m}(x;t) \exp(-2V(x;t)) \, {\rm d}x &= 0, \qquad \text{for\,}m\not\in \{n, n \pm 1\}, \\
    \int Q_{2n+1}(x;t)\,Q_{2m+1}(x;t) \exp(-2V(x;t)) \, {\rm d}x &= 0, \qquad \text{for\,}m\not\in \{n, n \pm 1, n\pm 2\}.
   \end{split}
\end{equation}

Furthermore, recall that according to \eqref{eq:p_ortho_laguerre}, each even-degree polynomial $P_{2n}(x;t)$ is mapped into the $n$-th degree polynomial $\widehat{P}_n(z;t)$, with $\widehat{P}(z;t)$  the family of monic polynomials orthogonal with respect to {\color{black}the Laguerre weight} $w_{-\frac{1}{2}}(z;t)$ (see \eqref{eq:semi-Laguerre_general}) by setting $x^2 = z$. Similarly by \eqref{eq:Ptildedef}  each odd-degree polynomial $P_{2n+1}(x;t)$ can be mapped into the $n$-th degree polynomial $\widetilde{P}_n(z;t)$, with $\widetilde{P}(z;t)$ the family of monic polynomials orthogonal with respect to {\color{black}the Laguerre weight} $w_{\frac{1}{2}}(z;t)$.
We can then introduce a similar mapping for skew-orthogonal polynomials defining
\begin{subequations}\label{eq:skew_Laguerre_poly}
\begin{align}
    \widehat{Q}_{n}(z;t)&:= Q_{2n}(x;t) \label{eq:skew_Laguerre_even}
    ,\\
   \widetilde{Q}_n(z;t) &:= \frac{Q_{2n+1}(x;t)}{x},
   \label{eq:skew_Laguerre_odd}
\end{align}
\end{subequations}
for any $n\geq 0$.
From Theorem \ref{thm:skew_quasi} we therefore have the following

\begin{corollary}\label{cor:quasi_Laguerre}
The monic polynomials $\widehat{Q}(z;t) = (\widehat{Q}_0(z;t), \widehat{Q}_1(z;t),\dots)^\top$ are quasi-orthogonal of order~$(1,1)$ with respect to the semi-classical Laguerre weight $w_{-\frac{1}{2}}(z;t)$. In particular, they can be expressed as 
\begin{equation}\label{eq:quasihat}
    \widehat{Q}_n = \widehat{P}_n + \xi_{n-1}\,\widehat{P}_{n-1}\,\qquad \forall ~n \in \mathbb{N}_0, 
\end{equation}
with polynomials $\widehat{P}(z;,t)$ as per \eqref{eq:p_ortho_laguerre} and $\{\xi_j(t)\}_{j\in\mathbb{N}_0}$ as per \eqref{eq:xi}.
Analogously the monic polynomials $\widetilde{Q}(z;t) = (\widetilde{Q}_0(z;t), \widetilde{Q}_1(z;t), \dots)^\top$ are quasi-orthogonal of order $(2,1)$ with respect to the semi-classical Laguerre weight $w_{\frac{1}{2}}(z;t)$. Specifically, we have 
\begin{equation}\label{eq:quasitilde}
    \widetilde{Q}_n = \widetilde{P}_n + \psi_{n-1}\,\widetilde{P}_{n-1}+\zeta_{n-2}\widetilde{P}_{n-2}\,\qquad \forall ~n \in \mathbb{N}_0, 
\end{equation}
with the understanding that all quantities indexed with negative indices {\color{black}in equations \eqref{eq:quasihat} and \eqref{eq:quasitilde}} are null. In the expression above, polynomials $\widetilde{P}(z; t)$ are defined in \eqref{eq:Ptildedef}, $\{\psi_j(t)\}_{j\in\mathbb{N}_0}$ in \eqref{eq:psi} and $\{\zeta_j(t)\}_{j\in\mathbb{N}_0}$ in \eqref{eq:zeta}.
\end{corollary}
We will further comment on the structure of skew-orthogonal polynomials as quasi-orthogonal ones in the next section.

\begin{remark}
    Given the definitions~\eqref{eq:skew_Laguerre_poly}, we look for a suitable expression of a skew-symmetric inner product for $\widehat{Q}_i(u;t)$ and $\widetilde{Q}_j(u;t)$, such that 
\begin{equation}
    \langle\,Q_{2i}(x;t), \,Q_{2j+1}(y;t)\,\rangle_{t} = \delta_{i,j} =:  \langle\,\widehat{Q}_{i}(u;t)\,|\,\widetilde{Q}_{j}(v;t)\,\rangle_{\ell} , \qquad \ell = 1, 2,3,4. 
\end{equation}
The choice to write the product 
$\langle~\cdot~|~\cdot~\rangle_{\ell}$ is not unique. In particular, it can be expressed either as the product $\langle~\cdot~,~\cdot~\rangle_t$ in terms of the same weight function for both sides, or in terms of two different weights, with the appropriate kernel $k_{\ell}(u,v)$. The fact that the map between $\{Q_{2i+1}(x;t)\}_{i\in\mathbb{N}_0}$ and $\widetilde{Q}(z;t)$ and the one from $\{Q_{2i}(x;t)\}_{\color{black}i\in\mathbb{N}_0}$ to $\widehat{Q}(z;t)$ in \eqref{eq:skew_Laguerre_poly} have a different functional dependency on $x$ is reflected by the form of the kernels $k_\ell$, $\ell =1,2,3,4$, which do not show symmetry or skew-symmetry properties  (i.e.\ $k_\ell(u,v)\neq \pm k_\ell(v,u)$).
We consider the change of variables $x^2=u$, $y^2=v$. One possibility is to write the product in terms of the  same weight for both variables $u$ and $v$ (Laguerre with $\lambda=\pm \frac{1}{2}$) as
\begin{subequations}
\begin{align}
\label{eq:skew_hat_tilde_1}
    \langle\,\widehat{Q}_{i}(u;t)\,|\,\widetilde{Q}_{j}(v;t)\,\rangle_{1} &= \frac{1}{2} \int_{\mathbb{R}_+^2} \widehat{Q}_{i}(u;t)\,\widetilde{Q}_{j}(v;t)\,w_{-\frac{1}{2}}\!(u;t)\,w_{-\frac{1}{2}}\!(v;t)\,k_1(u,v)\,{\rm d}u\,{\rm d}v,  \\
    \label{eq:skew_hat_tilde_2}
    \langle\,\widehat{Q}_{i}(u;t)\,|\,\widetilde{Q}_{j}(v;t)\,\rangle_{2} &= \frac{1}{2} \int_{\mathbb{R}_+^2} \widehat{Q}_{i}(u;t)\,\widetilde{Q}_{j}(v;t)\,w_{\frac{1}{2}}\!(u;t)\,w_{\frac{1}{2}}\!(v;t)\,k_2(u,v)\,{\rm d}u\,{\rm d}v,  
\end{align}
\end{subequations}
with the kernels $k_1(u,v)$ and $k_2(u,v)$ being respectively
    \begin{equation}\label{eq:ker12}
        k_1(u,v) = \frac{\sqrt{v}}{4}\,\text{\color{black}sgn}(v-u), \qquad k_2(u,v) = \frac{\text{\color{black}sgn}(v-u)}{4u\sqrt{v}}. 
    \end{equation}
We could adopt a more general form for the product  $\langle~\cdot~|~\cdot~\rangle_{\ell}$ involving two different weights. We have the following options
\begin{subequations}
\begin{align}\label{eq:skew_hat_tilde_3}
        \langle\,\widehat{Q}_{i}(u;t)\,|\,\widetilde{Q}_{j}(v;t)\,\rangle_{3} &= \frac{1}{2} \int_{\mathbb{R}_+^2} \widehat{Q}_{i}(u;t)\,\widetilde{Q}_{j}(v;t)\,w_{\frac{1}{2}}\!(u;t)\,w_{-\frac{1}{2}}\!(v;t)\,k_3(u,v)\,{\rm d}u\,{\rm d}v,  \\
        \label{eq:skew_hat_tilde_4}
        \langle\,\widehat{Q}_{i}(u;t)\,|\,\widetilde{Q}_{j}(v;t)\,\rangle_{4} &= \frac{1}{2} \int_{\mathbb{R}_+^2} \widehat{Q}_{i}(u;t)\,\widetilde{Q}_{j}(v;t)\,w_{-\frac{1}{2}}\!(u;t)\,w_{\frac{1}{2}}\!(v;t)\,k_4(u,v)\,{\rm d}u\,{\rm d}v,
    \end{align}
\end{subequations}
with the kernel $k_3(u,v)$ and $k_4(u,v)$ being 
\begin{equation}\label{eq:kuv_34}
        k_3(u,v) = \dfrac{\sqrt{v}}{4u}\,\text{\color{black}sgn}(v-u), \qquad k_4(u,v) = \dfrac{\text{\color{black}sgn}(v-u)}{4\sqrt{v}}.  
    \end{equation}
With the choice~\eqref{eq:skew_hat_tilde_3}
the polynomial $\widehat{Q}_i(u;t)$ ($\widetilde{Q}_i(u;t)$) appears as related to the Laguerre weight with $\lambda=\frac{1}{2}$ ($\lambda=-\frac{1}{2}$). However, $\widehat{Q}_i(u;t)$ ($\widetilde{Q}_i(u;t)$) is expressed as a combination of polynomials $\widehat{P}_i(u;t)$ ($\widetilde{P}_i(u;t)$) orthogonal with respect to Laguerre with $\lambda = -\frac{1}{2}$ ($\lambda =\frac{1}{2}$), as in~\eqref{eq:quasihat} (\eqref{eq:quasitilde}). This relation is restored in~\eqref{eq:skew_hat_tilde_4}. For each $\ell$, the product $\langle ~\cdot ~|~ \cdot ~\rangle_{\ell}$ behaves like a skew‑symmetric form when evaluated on pairs $(\widehat{Q}_i, \widetilde{Q}_j)$, i.e.\ it satisfies 
$\langle \,\widehat{Q}_i\,|\, \widetilde{Q}_j \,\rangle_{\ell} = -\langle \,\widetilde{Q}_j\,|\, \widehat{Q}_i \,\rangle_{\ell}$.
\end{remark}

\begin{remark}
It was shown in \cite{Adler2000_skew} that, if $w(x;t)$ in \eqref{eq:freud_weight}  is {\color{black} replaced by} the Gaussian weight $\exp(-\frac{x^2}{2})$ (with corresponding orthogonal polynomial $P(x;\,0)$ given by the standard Hermite polynomials), via a similar mapping to the one discussed above, skew-orthogonal polynomials can be expressed as 
\begin{subequations}
 \begin{align}
        Q_{2n}(x;\,0) &= P_{2n}(x;\,0), \\
        Q_{2n+1}(x;\,0) &= P_{2n+1}(x;\,0) - nP_{2n-1}(x;\,0) . \label{eq:oddQandPquadratic}
\end{align}
\end{subequations} 
for any $n\in \mathbb{N}_0$ and an appropriate shift \eqref{eq:Qodd_shift}. More insight on this scenario is provided in Appendix~\ref{sec:skewHermite}.

It can be checked with some algebra that if we consider a weight with a sextic potential of the form 
$$
w(x;t,\tau) = \exp\left(-V(x;t,\,\tau)\right),\qquad  V(x;t,\,\tau) = x^6 - \tau x^4 - t x^2,
$$ 
with polynomials $P(x;t,\tau) = (P_0(x;t,\tau), P_1(x;t,\tau),\dots)^{\top}$ orthogonal with respect to the symmetric weight $w(x;t,\tau)^2$ and  $Q(x;t,\tau) = (Q_0(x;t,\tau), Q_1(x;t,\tau),\dots)^{\top}$ skew-orthogonal with respect to $w(x;t,\tau)$, then the same type of mapping gives 
 \begin{align*}
         Q_{2n}(x;t,\tau) &= P_{2n}(x;t,\tau) + \phi_{n-1}^{(2)} \,P_{2n-2}(x;t,\tau)+ \phi_{n-1}^{(4)} \,P_{2n-4}(x;t,\tau), \\
         Q_{2n+1}(x;t,\tau) &= P_{2n+1}(x;t,\tau) + \phi_{n-1}^{(1)} \,P_{2n-1}(x;t,\tau) + \phi_{n-2}^{(3)}\,P_{2n-3}(x;t,\tau)+ \phi_{n-3}^{(5)}\,P_{2n-5}(x;t,\tau),
\end{align*}
with appropriate coefficients $\{\phi^{(\ell)}_j(t,\tau)\}_{j\in\mathbb{N}_0}$,  $\ell=1, \dots, 5$. 

We emphasise that orthogonal polynomials and the behaviour of $\{\beta_j(t,\tau)\}_{j\in\mathbb{N}_0}$ for the sextic Freud weight has been studied extensively, e.g.\ recently in \cite{Clarkson25} in relation to the discrete Painlev\'{e}-I hierarchy and in different regimes of the relevant sextic potential, and in~\cite{HermitianPRE} in the context of the Hermitian random matrix model. Skew-orthogonal polynomials and relevant recurrence relations will be object of future work.
These results can be generalised to $V(x)$ an even polynomial of higher degree. 
\end{remark}

\section{Recurrence relations for skew-orthogonal polynomials}
\label{sec:skew&quasi}
As shown by Theorem \ref{thm:skew_quasi} and Corollary \ref{cor:quasi_Laguerre}, monic skew-orthogonal polynomials with respect to the Freud weight $w(x;t)$~\eqref{eq:freud_weight} constitute two families of quasi-orthogonal polynomials~\cite{chihara,Draux01092016}, depending on their parity. In this section we therefore consider even and odd skew-orthogonal polynomials separately as two different families of quasi-orthogonal polynomials. Inspired by the approach discussed in \cite{IsmWan}, we exploit their quasi-orthogonal nature to formulate explicit recurrence relations. Indeed, in~\cite{IsmWan} the authors investigate classes of quasi-orthogonal polynomials defined in terms of orthogonal polynomials, providing a closed expression for the quasi-orthogonal polynomials, amongst other results. We follow the same strategy to determine a closed form for the skew-orthogonal polynomials with even indices $\{Q_{2j}(x;t)\}_{j\in\mathbb{N}_0}$ and odd indices $\{Q_{2j+1}(x;t)\}_{j\in\mathbb{N}_0}$ respectively. 

\subsection{Even skew-orthogonal polynomials}\label{sec:rec_even_skew}

We start by recalling equation \eqref{eq:two_step_three_terms}, obtained using the three-point relation in~\eqref{eq:3pts} twice, and evaluating it for orthogonal polynomials of even degrees $(2n-2)$ and $2n$:
\begin{subequations}
\begin{align}
 x^2\,P_{2n-2}&=P_{2n} + (\beta_{2n-2}+\beta_{2n-1})P_{2n-2} + \beta_{2n-3}\,\beta_{2n-2}\, P_{2n-4},       \label{eq:double3pts2n-2}\\
        x^2\,P_{2n}&=P_{2n+2} + (\beta_{2n}+\beta_{2n+1})P_{2n} + \beta_{2n-1}\,\beta_{2n}\, P_{2n-2}.
        \label{eq:double3pts2n}
        \end{align}
    \end{subequations}
Also recall that, according to~\eqref{eq:evenQandP}, we have 
\begin{subequations}\label{eq:QandPevenex}
\begin{align}
   Q_{2n+2} &= P_{2n+2} + \xi_{n} \,P_{2n},\label{eq:Q2np2}\\
    Q_{2n} &= P_{2n} + \xi_{n-1} P_{2n-2},\label{eq:Q2n}\\
    Q_{2n-2} &= P_{2n-2} + \xi_{n-2} \,P_{2n-4}.\label{eq:Q2nm2}
    \end{align}
    \end{subequations}
Moreover, \eqref{eq:double3pts2n-2} can be rearranged to provide us
with $P_{2n-4}(x;t)$ in terms of $P_{2n}(x;t)$ and $P_{2n-2}(x;t)$. 
Therefore, \eqref{eq:Q2n}-\eqref{eq:Q2nm2} take the matrix form {\color{black} $(Q_{2n},Q_{2n-2})^{\top}=A(x;t) (P_{2n},P_{2n-2})^{\top}$, or explicitly  }
\begin{equation}\label{eq:eventransform}
        \begin{pmatrix} Q_{2n} \\[1.5ex]
        Q_{2n-2}
        \end{pmatrix} = \begin{pmatrix} 
        1 & \xi_{n-1} \\
        -\dfrac{\xi_{n-2}}{\beta_{2n-3}\,\beta_{2n-2}}& 1+\dfrac{\xi_{n-2}}{\beta_{2n-3}\,\beta_{2n-2}}\left(x^2-(\beta_{2n-2}+\beta_{2n-1}) \right)
        \end{pmatrix} \,  \begin{pmatrix} P_{2n} \\[1.5ex]
        P_{2n-2}
        \end{pmatrix} \,{\color{black}.}
    \end{equation}
{\color{black} As the determinant {of \color{black}$A(x;t)$} is non zero, this matrix} can be inverted to obtain {\color{black} $(P_{2n},P_{2n-2})^{\top} = A^{-1}(x;t) (Q_{2n},Q_{2n-2})^{\top} $, or explicitly  }
{\color{black}
\begin{equation}\label{eq:Peven_Qeven}
    \begin{pmatrix} P_{2n} \\[2ex]
        P_{2n-2}
        \end{pmatrix} = \dfrac{1}{\det A} \begin{pmatrix} 
    1-\dfrac{\xi_{n-2} \left(\beta_{2 n-2}+\beta_{2 n-1}-x^2\right)}{\beta_{2 n-3} \beta_{2 n-2}} & -\xi_{n-1} \\[2ex]
     \dfrac{\xi_{n-2}}{\beta_{2 n-3} \beta_{2 n-2}} & 1 
        \end{pmatrix} \,  \begin{pmatrix} Q_{2n} \\[2ex]
        Q_{2n-2}
        \end{pmatrix}. 
\end{equation}
Furthermore, we note that equation~\eqref{eq:double3pts2n} can be solved for $P_{2n+2}(x;t)$ and considering the expression~\eqref{eq:Q2np2} for $Q_{2n+2}(x;t)$ we have 
\begin{equation}
\label{eq:Q2np2A}
    Q_{2n+2} = (x^2-\big(\beta_{2n}+\beta_{2n+1})+\xi_{n}\big)\,P_{2n}-\beta_{2n-1}\,\beta_{2n}\,P_{2n-2} .
\end{equation}    
Then \eqref{eq:Peven_Qeven} allows us to express both $P_{2n}$ and $P_{2n-2}$ in the right hand side of \eqref{eq:Q2np2A} as functions of $Q_{2n}$ and $Q_{2n-2}$ only, yielding the closed expression }
\begin{equation}
\begin{split} 
    &Q_{2n+2} = -\frac{\left( \beta_{2 n-3} \,\beta_{2 n-2} \,\beta_{2 n-1}\beta_{2 n}+\beta_{2 n-3}\, \beta_{2 n-2}\, \xi_{n-1} \left(x^2-\beta_{2 n}-\beta_{2 n+1}+\xi_{n}\right)\right)}{\beta_{2 n-3} \,\beta_{2 n-2}+\xi_{n-2} \,\left(x^2-\beta_{2 n-2}-\beta_{2 n-1}+\xi_{n-1}\right)}\,Q_{2n-2}\\
    &~~+\frac{\left(\beta_{2 n-3}\, \beta_{2 n-2}+\xi_{n-2} \left(x^2-\beta_{2 n-2}-\beta_{2 n-1}\right) \right)\!  \left(x^2-\beta_{2 n}-\beta_{2 n+1}+\xi_{n}\right)-\beta_{2 n-1}\, \beta_{2 n} \,\xi_{n-2}}{\beta_{2 n-3} \,\beta_{2 n-2}+\xi_{n-2} \,\left(x^2-\beta_{2 n-2}-\beta_{2 n-1}+\xi_{n-1}\right)}\,Q_{2n}.
    \label{eq:rec_Qeven_intermediate}
\end{split} 
\end{equation}
{\color{black} Within the expression above, we identify some structures which appear multiple times with different shifts in the variable $n$, which we denote as
\begin{equation}\label{eq:aux_funct_even}
    f_{n}(x;t):= \beta_{2n}(t)\,\beta_{2n-1}(t)+\xi_{n-1}(t)\,g_{n}(x;t),\qquad g_{n}(x;t) := x^2 - \beta_{2n}(t) - \beta_{2n+1}(t)+\xi_n(t) ,
\end{equation}
This observation allows us to rewrite equation \eqref{eq:rec_Qeven_intermediate} more compactly as detailed in the following statement:}
\begin{theorem}\label{thm:even_rec}
Let $\{Q_{2\ell}(x;t)\}_{\ell\in\mathbb{N}_0}$ be the subset of {\color{black} monic} polynomials of even degree of $Q(x;t)$ defined in~\eqref{eq:skew-ortho_monic} skew-{\color{black}orthogonal} with respect to the product~\eqref{eq:skew_product}. They satisfy the recursive equation
\begin{equation}\label{eq:even_closed_relation}
    f_{n-1}\,Q_{2n+2}= \big( f_{n-1}\,g_{n}-\xi_{n-2}\,f_{n} \big)\,Q_{2n} - \beta_{2n-3}\,\beta_{2n-2}\,f_{n}\, Q_{2n-2},  
\end{equation}
for $n\geq 2$, with $f_n = f_n(x;t)$ and $g_{n}=g_n(x;t)$ of the form~\eqref{eq:aux_funct_even}.
\end{theorem}

We emphasise that at this level the expressions we get are for even-indexed polynomials skew-orthogonal and orthogonal with respect to the (re-scaled) Freud weight in~\eqref{eq:skew_product} and \eqref{eq:product}, but these can be significantly simplified if formulated in terms of their semi-classical Laguerre counterpart.

Indeed, recall that according to \eqref{eq:p_ortho_laguerre}, by setting $x^2 = z$, each even-degree polynomial $P_{2n}(z;t)$ is mapped into the $n$-th degree polynomial $\widehat{P}_n(z;t)$, with $\widehat{P}(z;t)$ the family of monic polynomials orthogonal with respect to {\color{black}the Laguerre weight} $w_{-\frac{1}{2}}(z;t)$ (see \eqref{eq:semi-Laguerre_general}). Similarly, each even-degree polynomial $Q_{2n}(x;t)$ is mapped into the $n$-th degree polynomial $\widehat{Q}_n(z;t)$, with $\widehat{Q}(z;t)$ {\color{black} a} family of monic quasi-orthogonal polynomials of order (1,1) with respect to {\color{black}the Laguerre weight} $w_{-\frac{1}{2}}(z;t)$, as detailed in~\eqref{eq:skew_Laguerre_even} and Corollary~\ref{cor:quasi_Laguerre}. The following Corollary follows straightforwardly from Theorem \ref{thm:even_rec} and this mapping.
\begin{corollary}
Let $\widehat{Q}(z;t) = (\widehat{Q}_0(z;t), \widehat{Q}_1(z;t), \dots)^\top$ be the family of quasi-orthogonal polynomials with respect to the semi-classical Laguerre weight $w_{-\frac{1}{2}}(z;t)$ as per \eqref{eq:skew_Laguerre_even} and Corollary \eqref{cor:quasi_Laguerre}. They satisfy the following three-point recursive relation
\begin{equation}
    \widehat{f}_{n-1}\,\widehat{Q}_{n+1}= \big( \widehat{f}_{n-1}\,\widehat{g}_{n}-\xi_{n-2}\,\widehat{f}_{n} \big)\,\widehat{Q}_{n} - \widehat{b}_{n-1}\,\widehat{f}_{n}\, \widehat{Q}_{n-1},  
    \label{eq:rec_quasiortho_-1/2}
\end{equation}
for $n\geq 2$, with $\widehat{f}_n = \widehat{f}_n(z;t)$ and $\widehat{g}_n=\widehat{g}_n(z;t)$ defined as   
\begin{equation}
    \widehat{f}_n(z;t) = \widehat{b}_n(t)+\xi_{n-1}(t)\,\widehat{g}_{n}(z;t) , \qquad \widehat{g}_n(z;t) = z - \widehat{a}_n(t) + \xi_n(t) .
\end{equation}
Here, $\widehat{a}_n(t)$ and $\widehat{b}_n(t)$ have been defined in equation \eqref{eq:laguerre_coefficients} {\color{black} as, respectively, a sum and a product of coefficients $\beta_j(t)$ for suitable choices of $j$}. 
As usual, in \eqref{eq:rec_quasiortho_-1/2} we use the convention that $\widehat{Q}_j(z;t) =0$ for~$j<0$.
\label{cor:rec_Laguerre_even}
\end{corollary}
We note that the formula above is consistent with the general recursion relation for quasi-orthogonal polynomials of type \eqref{eq:quasiortho1st} as reported in \eqref{eq:quasiortho-rec}, keeping into account the recursion relation \eqref{eq:rec_Phat} satisfied by the polynomials $\widehat{P}(z;t)$.

\subsection{Odd skew-orthogonal polynomials}
\label{sec:odd_skew}
In this section, we proceed similarly to the previous one and focus on skew-orthogonal polynomials of odd degree. As previously mentioned, the approach is inspired by \cite{IsmWan}, with the main difference that in this case each odd-degree polynomial $Q_{2n+1}(x;t)$ is a combination of three orthogonal polynomials instead of two, which will lead to a slight generalisation of the method detailed there. 

Firstly, we recall equation
\eqref{eq:x2odd}, involving odd orthogonal polynomials only:
\begin{equation*}
    x^2 P_{2j+1} = P_{2j+3} + (\beta_{2j+2} + \beta_{2j+1})P_{2j+1} + \beta_{2j+1}\,\beta_{2j}\,P_{2j-1} ,
\end{equation*}
which can be evaluated for $j=\{n,n-1,n-2\}$ to find respectively
\begin{subequations}\label{eq:oddP2}
\begin{align}
\label{eq:oddP2piu1}
x^2 P_{2n+1} & = P_{2n+3} + (\beta_{2n+2} + \beta_{2n+1})P_{2n+1} + \beta_{2n+1}\,\beta_{2n}\,P_{2n-1},
\\
\label{eq:oddP2meno1}
    x^2 P_{2n-1}& = P_{2n+1} + (\beta_{2n} + \beta_{2n-1})P_{2n-1} + \beta_{2n-1}\,\beta_{2n-2}\,P_{2n-3},\\
\label{eq:oddP2meno3}
    x^2 P_{2n-3}& = P_{2n-1} + (\beta_{2n-2} + \beta_{2n-3})P_{2n-3} + \beta_{2n-3}\,\beta_{2n-4}\,P_{2n-5}. 
\end{align}
\end{subequations}
Moreover, notice that shifting equation \eqref{eq:oddQandP} we obtain 
\begin{subequations}
\label{eq:oddQandPshift}
 \begin{align}
 Q_{2n+3} &= P_{2n+3} + \psi_{n} \,P_{2n+1} + \zeta_{n-1}\,P_{2n-1},\label{eq:Q2np3P}\\
  Q_{2n+1} &= P_{2n+1} + \psi_{n-1} \,P_{2n-1} + \zeta_{n-2}\,P_{2n-3},\label{eq:Q2np1P}\\
 Q_{2n-1} &= P_{2n-1} + \psi_{n-2} \,P_{2n-3} + \zeta_{n-3}\,P_{2n-5}. \label{eq:Q2nm1P}
 \end{align}
 \end{subequations}
{\color{black}In~\eqref{eq:Q2np3P}} $P_{2n+3}$ can be expressed as a function of $P_{2n+1}$ {\color{black}and} $P_{2n-1}$ via~\eqref{eq:oddP2piu1}, and {\color{black} in~\eqref{eq:Q2nm1P}} $P_{2n-5}$ in terms of $P_{2n-3}$ {\color{black}and} $P_{2n-1}$ exploiting~\eqref{eq:oddP2meno3}. In this way, we determine {\color{black}the} closed system of equations
\begin{equation}
    \begin{pmatrix}
        Q_{2n+3} \\[.5ex] 
        Q_{2n+1} \\[.5ex] 
        Q_{2n-1} 
    \end{pmatrix} = {\color{black}B(x;t)} \begin{pmatrix}
        P_{2n+1} \\[.5ex]
        P_{2n-1} \\[.5ex]
        P_{2n-3}
    \end{pmatrix} 
\label{eq:PQoddMatrixMarta}
\end{equation}
where the matrix ${\color{black}B(x;t)}$ takes the form:
\begin{equation*}
    \hspace*{-1.75ex}{\color{black}B} = \begin{pmatrix}
        \psi_n+(x^2-\beta_{2n+2}-\beta_{2n+1})& \zeta_{n-1}-\beta_{2n+1}\,\beta_{2n} & 0 \\[1.5ex]
        1 & \psi_{n-1} & \zeta_{n-2} \\[1.5ex]
        0 & 1-\dfrac{\zeta_{n-3}}{\beta_{2n-3}\,\beta_{2n-4}} & \psi_{n-2}+\dfrac{\zeta_{n-3}(x^2-\beta_{2n-2}-\beta_{2n-3})}{\beta_{2n-3}\,\beta_{2n-4}} 
    \end{pmatrix}. 
\end{equation*}
{\color{black} As its determinant is non zero, the matrix {\color{black}$B$} can be inverted to} obtain the linear system of equations for $(P_{2n+1},P_{2n-1},P_{2n-3})^{\top}$ in terms of $(Q_{2n+3},Q_{2n+1},Q_{2n-1})^{\top}$. In particular, 
\begin{equation*}
    {\color{black}B^{-1}}=\frac{1}{\det {\color{black}B}} \begin{pmatrix}
        {\color{black}B}_{1,1} & {\color{black}B}_{1,2} & \zeta_{n-2} (\zeta_{n-1}-\beta_{2 n} \,\beta_{2n+1} ) \\[1.5ex]
        {\color{black}B}_{2,1} & {\color{black}B}_{2,2} & -\zeta_{n-2} (\psi_{n}+x^2 - \beta_{2n+1} - \beta_{2n+2} ) \\[1.5ex]
        1 - \dfrac{\zeta_{n-3}}{\beta_{2n-3}\, \beta_{2n-4}} & {\color{black}B}_{3,2} & {\color{black}B_{3,3}} \\
    \end{pmatrix},
\end{equation*}
with elements ${\color{black}B}_{i,j}$ given by
\begin{align*}
    {\color{black}B}_{1,1} &= \psi_{n-2} \,\psi_{n-1}-\zeta_{n-2}  + \frac{\zeta_{n-3} \Big(\zeta_{n-2} +\psi_{n-1} (x^2 - \beta_{2n-2} - \beta_{2n-3}) \Big)}{\beta_{2n-3} \,\beta_{2n-4}},\\
       {\color{black}B}_{1,2}&= -(\zeta_{n-1}-\beta_{2n} \,\beta_{2n+1} ) \Big(\psi_{n-2}+ \zeta_{n-3}\,\dfrac{x^2-\beta_{2n-2}-\beta_{2n-3}}{\beta_{2n-3}\,\beta_{2n-4}}\Big), \\
          {\color{black}B}_{2,1} &= - \Big(\psi_{n-2}+ \zeta_{n-3}\,\dfrac{x^2-\beta_{2n-2}-\beta_{2n-3}}{\beta_{2n-3}\,\beta_{2n-4}}\Big), \\ 
    {\color{black}B}_{2,2} &= \Big(\psi_{n-2}+\zeta_{n-3}\,\dfrac{x^2-\beta_{2n-2}-\beta_{2n-3}}{\beta_{2n-3}\,\beta_{2n-4}} \Big) (\psi_{n}+x^2 - \beta_{2n+1} - \beta_{2n+2} ), \\
     {\color{black}B}_{3,2} &= \Big(\frac{\zeta_{n-3}}{\beta_{2n-4}\, \beta_{2n-3}}-1\Big) (\psi_{n}+x^2 - \beta_{2n+1} - \beta_{2n+2} ), \\
     {\color{black}B_{3,3}}&= {\color{black}\beta_{2 n}\, \beta_{2n+1} - \zeta_{n-1} + \psi_{n-1} (\psi_{n}+x^2 - \beta_{2n+1} - \beta_{2n+2} ).}
\end{align*}

In order to determine the closed expression for odd skew-orthogonal polynomials, we consider first the next equation {\color{black} of the form}~\eqref{eq:oddP2}
    \begin{align}\label{eq:oddP2piu3}
    x^2\,P_{2n+3} &= P_{2n+5}+(\beta_{2n+4}+\beta_{2n+3})\,P_{2n+3}+\beta_{2n+3}\,\beta_{2n+2}\,P_{2n+1},      
\end{align}
and then the next equation {\color{black} of the form}~\eqref{eq:oddQandPshift}, i.e.
\begin{align}\label{eq:oddQ2piu5}
    Q_{2n+5} &= P_{2n+5}+\psi_{n+1}\,P_{2n+3}+\zeta_{n}\,P_{2n+1}.
\end{align}
From~\eqref{eq:oddP2piu3} we derive $P_{2n+5}$ in terms of $P_{2n+3}$ {\color{black}and} $P_{2n+1}$, so that in~\eqref{eq:oddQ2piu5} the right hand side is completely determined by $P_{2n+3}$ {\color{black}and} $P_{2n+1}$. From here, we get two different expressions for~\eqref{eq:oddQ2piu5}{\color{black}.
Exploiting ~\eqref{eq:oddP2piu1} we can express $P_{2n+3}$ as an appropriate combination of $P_{2n+1}$ and $P_{2n-1}$, yielding 
\begin{equation}\label{eq:Qfive1}
\begin{split} 
   Q_{2n+5}&= \bigl( x^4-x^2(\beta_{2n+3}+\beta_{2n+4}-\psi_{n+1})
   +\zeta_n -\beta_{2n+1}(x^2-\beta_{2n+3}-\beta_{2n+4}
   +\psi_{n+1}) \\
   &\qquad  - \beta_{2n+2}(x^2-\beta_{2n+4} + \psi_{n+1} ) \bigr) P_{2n+1}   - \beta_{2n}\,\beta_{2n+1} \bigl(x^2-\beta_{2n+3}-\beta_{2n+4}+\psi_{n+1}\bigr) P_{2n-1}
   \end{split}
\end{equation}
On the other hand, using a combination of~\eqref{eq:oddP2piu1} and~\eqref{eq:oddP2meno1} we can express $P_{2n+3}$ in terms of $P_{2n-1}$ and $P_{2n-3}$, which from \eqref{eq:oddP2piu1} implies, expressing $Q_{2n+5}$ in terms of $P_{2n+1}$ and $P_{2n-1}$,
\begin{equation}\label{eq:Qfive2}
\begin{split} 
   Q_{2n+5}&=  \Bigl(x^4 - x^2 (\beta_{2 n+3}+\beta_{2 n+4}- \psi_{n+1})+ \beta_{2 n} \beta_{2 n+1} (\beta_{2 n+3}+\beta_{2 n+4}-\psi_{n+1}-x^2)\\
   &~~
   +(x^2-\beta_{2 n}-\beta_{2 n-1}) \bigl(\zeta_{n}+\beta_{2 n+1} (\beta_{2 n+3}+\beta_{2 n+4}-\psi_{n+1}-x^2) \\
   &~~ +\beta_{2 n+2} (\beta_{2 n+4}-\psi_{n+1}-x^2)\bigr) \Bigr) P_{2 n-1} -  \beta_{2 n-2} \beta_{2 n-1}\Bigl( x^4 -x^2 (\beta_{2 n+3}+\beta_{2 n+4} - \psi_{n+1})  \\
   &~~
   +\zeta_{n}
   +\beta_{2 n+1} (\beta_{2 n+3}+\beta_{2 n+4}-\psi_{n+1}-x^2)+\beta_{2 n+2} (\beta_{2 n+4}-\psi_{n+1}-x^2)\Bigr)P_{2 n-3}. 
    \end{split} 
\end{equation}} 
Equating {\color{black}the right hand sides of \eqref{eq:Qfive1} and \eqref{eq:Qfive2}}, and writing the elements in $(P_{2n+1},P_{2n-1},P_{2n-3})^{\top}$ in terms of $(Q_{2n+3},\,Q_{2n+1},\,Q_{2n-1})^{\top}$ via the inverse matrix {\color{black} $B^{-1}$, we obtain a closed expression in $Q_{2n+3}$, $Q_{2n+1}$, and $Q_{2n-1}$ only. The resulting expression is rather complicated, however it can be written in terms of the following auxiliary functions appearing with different shifts in the variable $n$:
\begin{equation}
\label{eq:coeff_closed_skew_odd}
\begin{aligned}
c_n(t) &:= \zeta_{n-1}(t)-\beta_{2n+1}(t)\,\beta_{2n}(t),\qquad  d_n(t) := -c_{n-1}(t)\,c_{n}(t),\\
f_n(x;t) &:= -(\psi_{n}(t)+x^2-\beta_{2n+1}(t)-\beta_{2n+2}(t)),\\
g_n(x;t) &:= \psi_{n-2}(t)\,\beta_{2n-4}(t)\,\beta_{2n-3}(t)\,+\zeta_{n-3}(t)\left(x^2 - \beta_{2n-2}(t)-\beta_{2n-3}(t)\right).
\end{aligned}
\end{equation} 
This simplifies the closed expression, proving the following} 
\begin{theorem}\label{thm:odd_rec}
Let $\{Q_{2\ell+1}(x;t)\}_{\ell\in\mathbb{N}_0}$ be the subset of {\color{black} monic} polynomials of odd degree of $Q(x;t)$ defined in~\eqref{eq:skew-ortho_monic} skew-{\color{black}orthogonal} with respect to the product~\eqref{eq:skew_product}.
They satisfy the recursive equation
\begin{equation}
\label{eq:odd_closed_relation}
\begin{split}
&\big(d_{n}+f_{n-1} \,g_n \big) Q_{2n+3 } + \big(c_n\,g_n +c_{n-2}\,f_n \,\beta_{2n-1}\,\beta_{2n-2} +f_n\,g_n \left(f_{n-1}+\psi_{n-1}\right)\big)Q_{2n+1}\\
 &+\beta_{2n-4}\,\beta_{2n-3}\,\big(d_{n+1} + f_n\, g_{n+1}\big)Q_{2n-1} = 0,
\end{split}
\end{equation}
for $n\geq 3$, with $c_n = c_n(t)$, $d_n = d_n(t)$, $f_n = f_n(x;t)$ and $g_n = g_n(x;t)$ {\color{black}in \eqref{eq:coeff_closed_skew_odd}.} 
\end{theorem}
We point out that  relation~\eqref{eq:odd_closed_relation} is {\color{black} the counterpart of} \eqref{eq:even_closed_relation} for odd{\color{black}-degree} skew-orthogonal polynomials.
Incidentally, recall that, by setting $x^2 = z$, each odd-degree polynomials $P_{2n+1}(x;t)$ is mapped into $\widetilde{P}_n(z;t)$, with $\widetilde{P}(z;t)$ the family of monic polynomials orthogonal with respect to {\color{black}the Laguerre weight} $w_\frac{1}{2}(z;t)$ in \eqref{eq:semi-Laguerre_general}. Analogously, each odd-degree polynomial $Q_{2n+1}(x;t)$ is mapped into $\widetilde{Q}_n(z;t)$, with $\widetilde{Q}(z;t)$ the family of monic polynomials quasi-orthogonal of order (2,1) with respect to the {\color{black}Laguerre} weight $w_\frac{1}{2}(z;t)$, as per \eqref{eq:skew_Laguerre_odd} and Corollary \ref{cor:quasi_Laguerre}. Therefore, under this mapping, from Theorem \ref{thm:odd_rec} we have the following
\begin{corollary}\label{cor:rec_Laguerre_odd}
Let $\widetilde{Q}(z;t) =(\widetilde{Q}_0(z;t), \widetilde{Q}_1(z;t),\dots)^\top$ be the family of quasi-orthogonal polynomials with respect to the Laguerre weight $w_{\frac{1}{2}}(z;t)$ as per \eqref{eq:skew_Laguerre_odd} and Corollary \eqref{cor:quasi_Laguerre}. They satisfy the following three-point recursive relation
\begin{equation}
\label{eq:odd_closed_relation_Laguerre}
\begin{split}
&\big(\widetilde{d}_{n}+\widetilde{f}_{n-1}\, \widetilde{g}_n  \big) \widetilde{Q}_{n+1} + \Big(\widetilde{c}_n\,\widetilde{g}_n+\widetilde{c}_{n-2}\,\widetilde{f}_n \,\widetilde{b}_{n-1}+ \widetilde{f}_n\,\widetilde{g}_n \big(\widetilde{f}_{n-1}+\psi_{n-1}) \Big)\widetilde{Q}_{n}\\
 &+\widetilde{b}_{n-2}\big(\widetilde{d}_{n+1} + \widetilde{f}_n \,\widetilde{g}_{n+1}\big)\widetilde{Q}_{n-1} = 0,
\end{split}
\end{equation}
for $n\geq 3$, with $\widetilde{c}_n = \widetilde{c}_n(t)$, $\widetilde{d}_n = \widetilde{d}_n(t)$, $\widetilde{f}_n = \widetilde{f}_n(z;t)$ and $\widetilde{g}_n = \widetilde{g}_n(z;t)$  defined as 
\begin{equation}
\begin{aligned}
\widetilde{c}_n(t) &:= \zeta_{n-1}(t)-\widetilde{b}_{n}(t),\qquad\widetilde{d}_n(t) := -\widetilde{c}_{n-1}(t)\widetilde{c}_{n}(t),\\
\widetilde{f}_n(z;t) &:= -z + \widetilde{a}_{n}(t) - \psi_{n}(t), \\
\widetilde{g}_n(z;t) &:= \widetilde{b}_{n-2}(t)\,\psi_{n-2}(t)+\zeta_{n-3}(t)\left(z - \widetilde{a}_{n-1}(t)\right).
\end{aligned}
\end{equation}
\end{corollary}

\section{Conclusions and open questions}
Within this work, we have shown that polynomials {\color{black}skew-orthogonal} with respect to a quartic Freud weight can be expressed as suitable linear combinations of the relevant orthogonal ones (Theorem \ref{thm:skew_quasi}), with coefficients satisfying novel coupled recurrence relations \eqref{eq:recurrence_psi&zeta} and \eqref{eq:recurrence_xi}. From this, a new link between skew-orthogonal and quasi-orthogonal polynomials follows, allowing us to formulate the first examples of closed relationships involving skew-orthogonal polynomials only (Theorems \ref{thm:even_rec} and \ref{thm:odd_rec}). Furthermore, we point out that similar results hold also in the simplified scenario where, instead of a quartic Freud weight, we consider a Gaussian weight $w_{\text{G}}(x) = \exp({\color{black}-\frac{x^2}{2}})$, corresponding to Hermite orthogonal polynomials. This is briefly explained in Appendix \ref{sec:skewHermite}.

These results open a number of open questions and research avenues, some of which we discuss here. Recent work has been devoted to higher order Freud weights \cite{ClaJorKel2016, ClarJordaLour2023, Clarkson25}, and related orthogonal polynomials -- further investigation will be needed to generalise the results in our paper to these weights and further semi-classical ones. Moreover, though this is not explicitly explored here, skew-orthogonal polynomials have a deep relationship with  matrix ensembles, notably with symmetric and symplectic ones, and their underlying integrable lattice structure, i.e.\ the so-called Pfaff lattice \cite{vanM_Pfaff_skew, vanM_Pfaff_tau_function ,adler2002, Adler1999ThePL, BDM1, BDM2}. In particular, the recurrence relations \eqref{eq:recurrence_psi&zeta} and \eqref{eq:recurrence_xi} are expected to be related to the explicit solution of the Pfaff lattice for the symmetric matrix ensemble with quartic potential, and this link will be explored in future work.
Following up on this, our {\color{black} results} forge a novel relationship between quasi-orthogonal polynomials \cite{IsmWan, Bracciali} and symmetric and symplectic matrix ensembles, which might provide us with another set of tools to investigate them. Lastly, it is well known that, for orthogonal polynomials, relevant recursion coefficients~$\{\beta_i\}_{i\in\mathbb{N}_0}$ satisfy appropriate discrete Painlev\'e equations \cite{FokasIts, Clarkson25}. We hypothesise that analogous results might hold for skew-orthogonal polynomials, and that the recurrence relation \eqref{eq:recurrence_xi} for the parameter~$\{\xi_i\}_{i\in\mathbb{N}_0}$ might constitute an example of discrete Painlev\'e equation \cite{GraRam2004} in higher dimensions.

\section*{Acknowledgments}
We would like to thank G.\ Filipuk for significant insights. We are grateful to  P.A.\ Clarkson, A.\ Doliwa,  A.\ Loureiro, A.\ Moro, A.\ Stokes and R.\ Willox for helpful discussions. {\color{black}We also thank the referees for their insightful comments and feedback, helping us improve the quality of our manuscript.} The authors also gratefully acknowledge the GNFM
– Gruppo Nazionale per la Fisica Matematica, INdAM (Istituto Nazionale di Alta Matematica), for
supporting activities that contributed to the research presented in this paper.
\\\\

\begin{appendices}

\section{\texorpdfstring{Explicit expressions for ${Q}_{2n}(x;t)$ and ${Q}_{2n+1}(x;t)$}{qevenodd}}\label{app:formskew}
In this Appendix we provide the explicit form of of the first few monic polynomials skew-orthogonal with respect to the quartic Freud weight. Their coefficients are formulated in terms of $\{\beta_{j}(t)\}_{j \in \mathbb{N}_0}$ and $\{\xi_{j}(t)\}_{j \in \mathbb{N}_0}$, given by~\eqref{eq:recurrence_beta} and~\eqref{eq:recurrence_xi} respectively, with the initial points {\color{black}$\beta_0$, $\beta_1$ as in Section~\ref{sec:orthointro}}, and $\xi_0$, $\xi_1$ determined in Appendix~\ref{app:initial}.

The first examples of polynomials of even degree are:
\begin{align*}
    Q_0(x;t) &= 1 , \qquad  
    Q_2(x;t) = x^{2} + \xi_0 - \beta_{1}   , \\
    Q_4(x;t) &= x^{4} +  \big( {\color{black}\xi_1} - \beta_1 - \beta_2-\beta_3 \big)x^2 +\beta_1 \big(\beta_3 - {\color{black}\xi_1} \big)  , \\
    Q_6(x;t) &= x^{6} + \bigl({\color{black}\xi_2} - \beta_1 - \beta_2-\beta_3 - \beta_4-\beta_5 \bigr)x^4   + \bigl( \beta_5\,\beta_3 + \beta_5\,\beta_2 + \beta_5\,\beta_1 + \beta_4\,\beta_2 + \beta_4\,\beta_1 \\
   &\qquad + \beta_3\,\beta_1 - {\color{black}\xi_2}\big( \beta_{1}+ \beta_{2}+ \beta_{3} \big)\!\bigr) x^2 +\beta_{1}\,\beta_{3}\big({\color{black}\xi_2}-  \beta_{5} \big),\\
    Q_8(x;t) &= x^{8} + \big({\color{black}\xi_3}-\beta_{1}- \beta_{2}- \beta_{3}- \beta_{4}- \beta_{5}- \beta_{6}- \beta_{7}\!\bigr)x^{6}     
    +
    \bigl( \beta_7\,\beta_5 + \beta_7\,\beta_4 + \beta_7\,\beta_3 + \beta_7\,\beta_2 \\
   &  + \beta_7\,\beta_1 
   + \beta_6\,\beta_4 + \beta_6\,\beta_3 + \beta_6\,\beta_2 + \beta_6\,\beta_1 
   + \beta_5\,\beta_3 + \beta_5\,\beta_2 + \beta_5\,\beta_1 + \beta_4\,\beta_2 + \beta_4\,\beta_1 + \beta_3\,\beta_1   \\
   &
   -{\color{black}\xi_3} \big(\beta_{1}+ \beta_{2}+ \beta_{3}+ \beta_{4}+ \beta_{5} \big)  \bigr) x^4  
   + \bigl( {\color{black}\xi_3} \big(\beta_5\,\beta_3 + \beta_5\,\beta_2 + \beta_5\,\beta_1 + \beta_4\,\beta_2 + \beta_4\,\beta_1 + \beta_3\,\beta_1 \big) \\
   &
   - \beta_7 \, \beta_5 \, \beta_3
   - \beta_7 \, \beta_5 \, \beta_2
   - \beta_7 \, \beta_5 \, \beta_1
   - \beta_7 \, \beta_4 \, \beta_2
   - \beta_7 \, \beta_4 \, \beta_1
   - \beta_7 \, \beta_3 \, \beta_1
   - \beta_6 \, \beta_4 \, \beta_2
   - \beta_6 \, \beta_4 \, \beta_1
   \\
   &
   - \beta_6 \, \beta_3 \, \beta_1
   - \beta_5 \, \beta_3 \, \beta_1  \bigr) x^2 +\beta_{1}\,\beta_{3}\,\beta_{5}\big(\beta_{7}- {\color{black}\xi_3} \big) ,
\end{align*} 
and for odd degree polynomials: 
\begin{align*}
    Q_1(x;t) &= x , \qquad 
    Q_3(x;t) = x^{3} + \Big( \frac{\xi_0\,\xi_1}{\beta_{1}} - \frac{{\color{black}\xi_0}}{4\beta_1\,\beta_2} - \beta_{1}- \beta_{2} \Big)x   , \\[.5ex]
    Q_5(x;t) &= x^{5}  + 
    \Big(  \frac{\xi_1\,\xi_2}{\beta_{3}} - \frac{2 {\color{black}\xi_1}}{4\beta_3\,\beta_4} - \beta_{1}- \beta_{2} - \beta_{3}- \beta_{4} \Big)x^3  \\
    &~~
    - \Big(  \Big(\frac{\xi_1\,\xi_2}{\beta_3}-\frac{2 {\color{black}\xi_1}}{4\beta_{3}\,\beta_4} \Big) (\beta_{1}+\beta_2 ) +\beta_4\,{\color{black}\xi_1} - \beta_4 \, \beta_2 - \beta_4 \beta_1 - \beta_3 \, \beta_1 \Big) x,        \\
    Q_7(x;t) &= x^{7} +  \Big( \frac{\xi_2\,\xi_3}{\beta_{5}} - \frac{3{\color{black}\xi_2}}{4\beta_5\,\beta_6} - \beta_{1}- \beta_{2}- \beta_{3}- \beta_{4}- \beta_{5}- \beta_{6} \Big)x^5 \\ 
     &- \Big(  \Big(\frac{\xi_2\,\xi_3}{\beta_{5}}-\frac{3{\color{black}\xi_2}}{4\beta_5\,\beta_6}\Big)\Big( \beta_1+\beta_2+\beta_3+\beta_4 \Big)  
     + \beta_8\,{\color{black}\xi_2}  
     -\beta_6\, \beta_4 -\beta_6\, \beta_3 -\beta_6\, \beta_2 \\
     & -\beta_6\, \beta_1 -\beta_5\, \beta_3         -\beta_5\, \beta_2 -\beta_5\, \beta_1 -\beta_4\, \beta_2 -\beta_4\, \beta_1 -\beta_3\, \beta_1  \Big)x^3 \\
     &+ \Big(  \Big(\frac{\xi_2\,\xi_3}{\beta_{5}}-\frac{3{\color{black}\xi_2}}{4\beta_5\,\beta_6}\Big) \Big( \beta_4\,\beta_2 + \beta_4\,\beta_1 + \beta_3\,\beta_1  \Big)      +\beta_8\,\,{\color{black}\xi_2} - \beta_6\,\beta_4\,\beta_2 - \beta_6\,\beta_4\,\beta_1 \\
     & - \beta_6\,\beta_3\,\beta_1  - \beta_5\,\beta_3\,\beta_1 \Big)x . 
\end{align*}

\section{\texorpdfstring{Initial values for the recurrence relation for $\xi_n$}{xin}}
\label{app:initial}
The recurrence relation \eqref{eq:recurrence_xi} plays a pivotal role in this work, as it allows us to determine the coefficients $\xi_n(t)$, underpinning the formulation of the skew-orthogonal polynomials as quasi-orthogonal ones. This is a three-point recurrence relation, and as such necessitates of two initial values as starting point. This Appendix is devoted to explicit calculations evidencing how to calculate $\xi_0(t)$ and $\xi_1(t)$. These can be implemented with appropriate mathematical software for any given value of the parameter $t$.
\subsection{\texorpdfstring{Evaluation of $\xi_0(t)$}{xi0}}\label{app:xi0}
Firstly, note that $Q_0(x;t) = P_0(x;t)=1$ and $Q_1 (x;t) = P_1(x;t)= x$. Moreover, recall that, based on \eqref{eq:evenQandP} we have
$$
Q_2(x;t) = P_2(x;t) + \xi_0(t) P_0(x;t) = P_2(x;t) + \xi_0(t).
$$
Therefore, we need to determine the explicit form of $P_2(x;t)$. Due to the symmetry of the scalar product \eqref{eq:product} and the recurrence relation \eqref{eq:3pts} we have
\begin{equation}
P_2(x;t) = x^2 -\beta_1(t),
\label{eq:P2}
\end{equation}
with $\{\beta_j(t)\}_{j\in\mathbb{N}_0}$ defined in \eqref{eq:alpha&beta} and \eqref{eq:beta_def}.
Thus $Q_2(x;t)$ takes the form $Q_2(x;t) = x^2 +\xi_0(t) -\beta_1(t)$. We impose that $\langle Q_2, Q_1 \rangle_t = 0$ to evaluate $\xi_0(t)$ and find
\begin{equation}
\langle x^2, y\rangle_t +(\xi_0(t)-\beta_1(t))\langle 1,y\rangle_t = 0~~{\color{black} \implies }~~ 
 \xi_0(t) = -\frac{\langle x^2, y\rangle_t }{\langle 1,y\rangle_t}+\beta_1(t).
\label{eq:xi0a}
\end{equation}
This can be evaluated for specific values of $t$ with standard mathematical software. 
Also note that \eqref{eq:xi0a} implies that
\begin{equation}
Q_2(x;t) = x^2 -\frac{\langle x^2,y\rangle_t}{\langle 1,y\rangle_t}.
\label{eq:Q2}
\end{equation}

We remark that there is a different, yet equivalent, expression for $\xi_0(t)$. Note that according to equation \eqref{eq:xi} we also have
\begin{equation}
\xi_0(t) = \frac{4r_1(t)}{h_0(t)} 
\label{eq:xi0b}
\end{equation}
{\color{black}with $h_0(t) = \left(1,1\right)_t$ and $r_1(t) = \langle Q_2, Q_3\rangle_t$}. To evaluate $r_1(t)$, the explicit form of $Q_3(x;t)$ is needed. By construction, $Q_3(x;t) = x^3 + \nu(t) x$, with $\nu(t)$ to be determined. We impose orthogonality with $Q_0(x;t)$ to find
$$
\langle 1, y^3\rangle_t + \nu(t)\langle 1,y\rangle_t = 0, 
$$
which implies  $\nu(t) =-\frac{\langle 1,y^3\rangle_t}{\langle 1,y\rangle_t}.$
Thus exploiting \eqref{eq:Q2} we have
{\color{black} 
\begin{equation}
  r_1(t) = \langle x^2, y^3\rangle_t  -\frac{\langle 1,y^3\rangle_t}{\langle 1,y\rangle_t}\langle x^2,y\rangle_t
\end{equation}}
which can be evaluated for a given choice of the parameter $t$. The two expressions \eqref{eq:xi0a} and \eqref{eq:xi0b} are equivalent and it is straightforward to confirm this by checking their values for different choices of $t$. 

\subsection{\texorpdfstring{Evaluation of $\xi_1(t)$}{xi}}
To evaluate $\xi_1(t)$, we exploit \eqref{eq:recxi_withr} evaluated for $\ell = 0$ which, recalling that  by definition $\mathcal{q}_{2,0}(t) = \xi_0(t)$ and $\mathcal{q}_{4,2}(t) = \xi_1(t)$, takes the form
\begin{equation}
    1 = \frac{r_0}{2 h_1} + \frac{r_0}{h_2}\xi_0 - 4 \frac{r_0}{h_1}\xi_1 \xi_0.
\end{equation}
This can be solved for $\xi_1(t)$ as
\begin{equation}
    \xi_1 = \frac{1}{4\xi_0}\left({\color{black}\frac{1}{2}}- \frac{h_1}{r_0}\right)  + \frac{1}{4 \beta_2}
\end{equation}
with $r_0 = \langle 1, y\rangle_t$ and $h_1 = \left(P_1,P_1\right)_t$. 
As $\xi_0(t)$ can be evaluated as detailed in Appendix \ref{app:xi0}, $\xi_1(t)$ can be evaluated for any given value of $t$ with appropriate mathematical software.

We also illustrate another way to compute $\xi_1(t)$ which does not rely on \eqref{eq:recxi_withr}.
Recall that $Q_4(x;t) = x^4 + \gamma_2(t) x^2 + \gamma_0(t)$, with $\gamma_2(t)$ and $\gamma_0(t)$ to be determined. Similarly, $Q_3(x;t) = x^3 + \nu(t)x$, with $\nu(t) = -\frac{\langle 1 ,y^3\rangle_t}{\langle 1,y\rangle_t }$, as shown in the previous section. To determine $\gamma_2(t)$ and $\gamma_0(t)$, we require that $\langle Q_4, Q_1\rangle_t = 0$, and $\langle Q_4,  Q_3\rangle_t = 0$.
This provides the following system of linear equations for $\gamma_0(t)$ and~$\gamma_2(t)$:
\begin{equation}
\begin{split}
    \langle x^4, y\rangle_t +\gamma_2 \langle x^2,y\rangle_t + \gamma_0 \langle 1,y\rangle_t &= 0\\
    \langle x^4, y^3\rangle_t + \nu \langle x^4, y\rangle_t + \gamma_2 \left(\langle x^2, y^3\rangle_t + \nu \langle x^2,y\rangle_t\right) + \gamma_0 \left(\langle 1,y^3\rangle_t + \nu \langle 1, y\rangle_t\right) &=0
    \end{split}
\end{equation}
Therefore we have
\begin{equation}
    \gamma_2 = \frac{\langle 1,y\rangle_t \langle x^4, y^3\rangle_t - \langle 1,x^3\rangle_t\langle x^4, y\rangle_t}{\langle 1, y^3\rangle \langle x^2, y\rangle_t - \langle 1,y\rangle_t \langle x^2, y^3 \rangle_t},\;\;\;\;\gamma_0 =  \frac{\langle x^2,y^3\rangle_t \langle x^4, y\rangle_t - \langle x^2,y\rangle_t\langle x^4, y^3\rangle_t}{\langle 1, y^3\rangle_t \langle x^2, y\rangle_t - \langle 1,y\rangle_t \langle x^2, y^3 \rangle_t},
\end{equation}
completely characterising $Q_4(x;t)$. 
Recall that 
\begin{equation}
Q_4(x;t) = P_4(x;t) + \xi_1(t) P_2(x;t).
\label{eq:Q4}
\end{equation}
 To find $\xi_1(t)$ is thus necessary to make $P_4(x;t)$ and $P_2(x;t)$ explicit. $P_2(x;t)$ is reported in \eqref{eq:P2}. Exploiting the recurrence relation \eqref{eq:3pts} we have
 \begin{equation}
     P_4(x;t) = x^4 - (\beta_1(t)+\beta_{\color{black}2}(t)+\beta_{\color{black}3
     }(t))x^2 + \beta_1(t)\,\beta_{3}(t).
     \label{eq:P4}
 \end{equation}
From \eqref{eq:Q4} it follows that 
\begin{equation}
    Q_4(x;t) = x^4 - \left(\beta_1(t)+\beta_{\color{black}2}(t)+\beta_{\color{black}3}-\xi_1(t)\right) x^2 + \beta_1(t)(\beta_3(t) -\xi_1(t)), 
\end{equation}
therefore
\begin{equation}
    \gamma_2 =  \xi_1- \beta_1-\beta_{\color{black}2}-\beta_{\color{black}3},\qquad \gamma_0 = \beta_1(\beta_3 -\xi_1).
\end{equation}

Solving these two equations for $\xi_1(t)$, we have found other two equivalent expressions to calculate its value,
\begin{equation}
    \xi_1 = \gamma_2+\beta_1 + \beta_{\color{black}2} + \beta_{\color{black}3}, \qquad \xi_1 = \beta_3 -\frac{\gamma_0}{\beta_1}.
\end{equation}
which can be easily evaluated with appropriate mathematical software.

\section{Recurrence relations for skew-Hermite polynomials}\label{sec:skewHermite}
 
We conclude the paper by stating the analogue of Theorems \ref{thm:even_rec} and \ref{thm:odd_rec} for the skew-orthogonal version of the Hermite polynomials, constructed for the Gaussian weight 
\begin{equation}\label{eq:gauss}
    w_{\text{G}}(x) = \exp\left(-\frac{x^2}{2}\right),
\end{equation}
and interpreted as quasi-orthogonal polynomials. In the following, we combine insights from~\cite{Adler2000_skew,IsmWan,BDM2}. 

Firstly, we introduce the skew-symmetric inner product with weight \eqref{eq:gauss}, defined in~\cite{BDM2} and that we denote with $\langle~\cdot~,~\cdot~\rangle_{\text{G}}$, i.e.
\begin{equation}
\langle\,f,\,g\,\rangle_{\text{G}} = \int_{\mathbb{R}^2} f(x)\,g(y)\,\text{\color{black}sgn}(y-x)\,w_{\text{G}}(x)\,w_{\text{G}}(y)\,\text{d}x\,\text{d}y
\label{eq:skew-prod_gauss}
\end{equation}
 We define the vector $S(x)=(S_{0}(x), S_1(x), \dots)^{\top}$ of monic polynomials skew-orthogonal with respect to $w_{\text{G}}(x)$ satisfying the following
\begin{equation*}
    \langle\,S_{2i},\,S_{2j+1}\,\rangle_{\text{G}} = - \langle\,S_{2j+1},\,S_{2i}\,\rangle_{\text{G}} = \delta_{i,j}\,k_j, \qquad \langle\,S_{2i},\,S_{2j}\,\rangle_{\text{G}} = 0 = \langle\,S_{2i+1},\,S_{2j+1}\,\rangle_{\text{G}},
\end{equation*}
with $\{k_j\}_{j\in\mathbb{N}_0}$ appropriate normalising coefficients. Furthermore, we introduce the vector of monic polynomials $O(x)=(O_{0}(x), O_1(x), \dots)^{\top}$ orthogonal with respect to the weight $w_{\text{G}}(x)^2 = \exp(-x^2)$, and the inner product $(~\cdot~,~\cdot~)_{\text{G}}$ such that: 
\begin{equation*}
    (\,f,\,g\,)_{\text{G}} = \int_{\mathbb{R}} f(x)\,g(x)\,w_{\text{G}}(x)^2\,\text{d}x,\qquad 
    (\,O_{i},\,O_{j}\,)_{\text{G}} =  \delta_{i,j}\,\ell_j,
\end{equation*}
with $\{\ell_j\}_{j\in\mathbb{N}_0}$ normalising coefficients. For any $n \in \mathbb{N}_0$ the polynomial $O_n(x)$ is proportional to the standard Hermite polynomial $H_n(x)=2^n\, O_n(x)$.
The three-point relation for the polynomials $O(x)$ is 
\begin{equation}\label{eq:threepoint_gauss}
   x\, O_n(x) = O_{n+1}(x)+\frac{n}{2}\,O_{n-1}(x).
\end{equation}  
As explicitly given in~\cite{Adler2000_skew}, the skew-orthogonal polynomials associated with the classical Gaussian weight are expressed as combinations of orthogonal ones as
\begin{subequations}
\begin{align}
\label{eq:gauss_pari}
    S_{2n}(x) &= O_{2n}(x) , \\ 
\label{eq:gauss_dispari}
    S_{2n+1}(x) &= O_{2n+1}(x) -n\,O_{2n-1}(x). 
\end{align}
\end{subequations}
Using twice~\eqref{eq:threepoint_gauss} on the subset of orthogonal polynomials with even index and of odd index, we obtain the following relations: 
\begin{subequations}
\begin{align}
       x^2\,O_{2n} &= O_{2n+2}+\left(2n+\frac{1}{2}\right)\,O_{2n}+n\left(n-\frac{1}{2}\right)\,O_{2n-2} , \\
        x^2\,O_{2n-1} &= O_{2n+1} + \left(2n-\frac{1}{2}\right)\,O_{2n-1}+(n-1)\left(n-\frac{1}{2}\right)\,O_{2n-3}. \label{eq:hermite_odd}
    \end{align}
\end{subequations}

From~\eqref{eq:gauss_pari}, the closed relation for the skew-orthogonal polynomials with even index is then simply 
\begin{equation}
    S_{2n+2} = \left(x^2 -2n-\frac{1}{2}\right)S_{2n}-n\left(n-\frac{1}{2}\right)S_{2n-2}.
\end{equation}
To determine the closed relation for the skew-orthogonal polynomials with odd index, we first build the matrix mapping from $(O_{2n+1},O_{2n-1})^{\top}$ into $(S_{2n+1},S_{2n-1})^{\top}$, exploiting \eqref{eq:gauss_dispari} and \eqref{eq:hermite_odd}
\begin{equation}\label{eq:hermite_transf}
    \begin{pmatrix}
        S_{2n+1} \\[1.75ex] 
        S_{2n-1}
    \end{pmatrix} = \begin{pmatrix}
        1 & -n \\ 
        \dfrac{1}{n-\frac{1}{2}} & \dfrac{-x^2+3n-1}{n-\frac{1}{2}} 
    \end{pmatrix} \begin{pmatrix}
        O_{2n+1} \\[1.75ex]
        O_{2n-1} 
    \end{pmatrix}.
\end{equation}
We then invert it to obtain the expressions for $(O_{2n+1},O_{2n-1})^{\top}$ in terms of $(S_{2n+1},S_{2n-1})^{\top}$.
Notice that, by shifting \eqref{eq:gauss_dispari} we have
\begin{equation}
    S_{2n+3}(x) = O_{2n+3}(x) -(n+1)\,O_{2n+1}(x).
\end{equation}
We can then exploit \eqref{eq:hermite_odd} shifted by one to express $O_{2n+3}(x)$ in terms of $O_{2n+1}(x)$ and $O_{2n-1}(x)$ to find
\begin{equation}
    S_{2n+3}(x) = \left(x^2 - \left(3n + \frac{5}{2}\right)\right)O_{2n+1}(x) -n\left(n+\frac{1}{2}\right)\,O_{2n-1}(x).
\end{equation}
Proceeding similarly to Section \ref{sec:rec_even_skew}, exploiting the inverse of transformation \eqref{eq:hermite_transf} with some algebra we find the following closed relation
\begin{equation}
\begin{split} 
    &\qquad  2(x^2 -4n + 1) S_{2n+3} = \\
    &=\big(2x^4 - 3(1 + 4n)x^2 +8n(2n+1) -5\big)S_{2n+1}-n(2n-1)(x^2 -4n-3)S_{2n-1}.
\end{split} 
\end{equation}

\end{appendices}

\bibliography{biblio}
\bibliographystyle{stylesort}

\end{document}